\documentclass[12pt]{article}
\usepackage{amsbsy}
\usepackage{amsfonts}
\usepackage{amsmath}
\usepackage{amssymb}
\usepackage{apacite}
\usepackage{setspace}
\usepackage{graphicx}
\usepackage{bm,multicol}
\usepackage[top=3cm, bottom=3cm, left=3cm, right=3cm]{geometry}
\usepackage{natbib}
\usepackage[T1]{fontenc}
\usepackage[utf8]{inputenc}
\usepackage{authblk}
\usepackage{xcolor}

\newtheorem{lemma}{Lemma}

\newtheorem{definition}{Definition}

\newtheorem{remark}{Remark}
\newtheorem{example}{Example}
\newcommand{\vet}[1]{{\ensuremath{\mbox{\boldmath $#1$}}}}
\doublespace
\newcommand{\newc}{\newcommand}
\newc{\N}{\mbox{N}}
\newc{\1}{\bf{1}}
\def\hatth{\widehat{\bm{\theta}}}
\def\th{\bm{\theta}}
\def\Om{\bm{\Omega}}
\def\J{\mathcal{I}_\theta^{-1}}
\def\I{\mathcal{I}_\theta}
\def\normth{\lVert \widehat{\bm{\theta}} \rVert}

\begin{document}
\title{On the Ubiquity of Information Inconsistency for Conjugate Priors}

\author{J. Mulder, J. O. Berger, V\'ictor Pe\~na and M. J. Bayarri}

\thispagestyle{empty} \maketitle

\begin{abstract}
Informally, `Information Inconsistency' is the property that has been observed
in many Bayesian hypothesis testing and model selection procedures whereby the
 Bayesian conclusion does not become definitive when the data seems to
 become definitive. An example is that, when performing a $t$-test using standard conjugate priors, the Bayes factor
  of the alternative hypothesis to the null hypothesis remains bounded as the $t$ statistic grows to infinity.
  This paper shows that information inconsistency is ubiquitous in Bayesian hypothesis testing under conjugate priors.
  Yet the title does not fully describe the paper, since we also show that theoretically recommended priors,
  including scale mixtures of conjugate priors and adaptive priors, are information consistent.
  Hence the paper is simply a forceful warning that use of conjugate priors in testing and model selection
  is highly problematical, and should be replaced by the information consistent alternatives.
\end{abstract}

\section{Introduction}
When testing a precise null hypothesis $H_0$ against an unrestricted alternative hypothesis $H_1$, a common Bayesian tool is the Bayes factor, $B_{10}$, which quantifies the relative evidence (or odds) from the data for $H_1$ against $H_0$. A Bayes factor is called {\em information inconsistent} if, when the evidence for the alternative hypothesis appears to be overwhelming (in the sense that the observed effect under
the alternative hypothesis becomes arbitrarily large), the Bayes factor converges to a constant $B^*<\infty$. This conflicting behavior is also referred to as the information paradox \citep{Liang:2008}.

\begin{example}
A typical example of an information inconsistent Bayes factor is when using Zellner's (1986) $g$-prior for testing the regression coefficients in a linear regression model $\textbf{y}=\gamma\textbf{1}_n+\textbf{X}_{\boldsymbol{\theta}}\boldsymbol{\theta}+\boldsymbol{\epsilon}$, with $\boldsymbol{\epsilon}\sim N(\textbf{0},\sigma^2\textbf{I}_n)$, where $\textbf{y}$ is a vector containing the $n$ responses, $\gamma$ is the intercept, $\textbf{X}_{\boldsymbol{\theta}}$ is a $n\times r_1$ matrix containing the explanatory variables, $\boldsymbol{\theta}$ is a vector with the $r_1$ unknown coefficients that are tested, $\sigma^2$ is the unknown error variance, $\textbf{1}_n$ is a vector of length $n$ with ones, and $\textbf{I}_n$ is the identity matrix of size $n$. When testing $H_0:\vet{\theta}=\vet{0}$ versus $H_1:\vet{\theta}\not =\vet{0}$ with the $g$-prior, $\pi_0(\gamma,\sigma^2)\propto\sigma^{-2}$ and $\pi_1(\boldsymbol{\theta} \mid \gamma,\sigma^2)=N_{\boldsymbol{\theta} \mid \gamma,\sigma^2}(\textbf{0},g \sigma^2 (\textbf{X}_{\boldsymbol{\theta}}'\textbf{X}_{\boldsymbol{\theta}})^{-1})$ and $\pi_1(\gamma,\sigma^2)\propto\sigma^{-2}$, for some fixed $g>0$, the Bayes factor goes to $(1+g)^{n-r_1-1}<\infty$ as the evidence against $H_0$ accumulates in the sense that $|\hat{\boldsymbol{\theta}}|\rightarrow\infty$ \citep[see also,][]{Berger:2001}. Furthermore, it has also been reported that the $g$-prior is information inconsistent when testing one-sided hypotheses \citep{Mulder:2014a}.
\end{example}

In comparison to large sample inconsistency, which occurs when the evidence for the true hypothesis against another hypothesis does not go to infinity as the sample size grows, information inconsistency has not received much attention in the literature. In our view, both types of inconsistency are undesirable and should be avoided in general testing procedures. The goal of this paper is therefore to explore information inconsistency in a more general setting. We will consider improper as well as proper priors; conjugate priors, scale mixtures of conjugate priors, independent priors, and adaptive priors; and precise null hypothesis testing, one-sided hypothesis testing, and multiple hypothesis testing. Throughout the paper we also consider variations of Zellner's $g$ prior (e.g., fixed $g$ priors, mixtures of $g$ priors, and adaptive (data-based) $g$ priors) as this class of priors is commonly observed in the literature. One of the main questions we want to address is whether information consistent Bayes factors can be obtained using `standard' conjugate or independent semi-conjugate priors, or whether more sophisticated scale mixtures or adaptive priors are needed. We also explore whether there are any practical consequences by investigating when information inconsistency starts to manifest itself and what the limiting value of the Bayes factor is.

The paper is organized as follows. First the linear regression model with dependent errors as well as some notation are introduced (Section \ref{section2}). Subsequently Section \ref{nullHypo} explores information consistency when testing a precise hypothesis using various prior specifications, followed by one-sided hypothesis tests in Section \ref{nonnestedHypo}, and a multiple hypothesis test in Section \ref{multipleHypo}. We end the paper with some conclusions and recommendations in Section \ref{concl}.

\section{The Linear Regression Model with Dependent Errors}\label{section2}
Throughout this paper the focus shall be on the linear regression model with dependent errors,
\begin{equation}
\textbf{y}=\textbf{X}_{\boldsymbol{\beta}}\boldsymbol{\beta}+\vet{\epsilon}, \mbox{ with }\vet{\epsilon}\sim N(\textbf{0},\sigma^2\vet{\Sigma}),
\label{model1}
\end{equation}
where the vector $\textbf{y}$ of length $n$ contains the responses, $\textbf{X}_{\boldsymbol{\beta}}=[\textbf{x}_1~\ldots~\textbf{x}_K]$ is an $n\times K$ matrix containing the $K$ predictor variables which are regressed on the $K$ unknown regression coefficients in $\vet{\beta}$ ($n > K$), $\vet{\epsilon}$ is a normally distributed error vector, $\sigma^2$ is an unknown common variance, and $\vet{\Sigma}$ is a known covariance matrix.

Three different types of hypothesis tests will be considered. First, we consider the classical null hypothesis test of a set of linear restrictions on $\bm{\beta}$ against an unrestricted alternative, i.e., $H_0:\textbf{R}\boldsymbol{\beta}=\textbf{0}$ versus $H_1:\textbf{R}\boldsymbol{\beta}\not =\textbf{0}$,
where $\textbf{R}$ is an $r_1\times K$ matrix with known constants ($r_1\le K$). Second, we consider the equivalent one-sided hypothesis test of
$H_0:\textbf{R}\boldsymbol{\beta}\le\textbf{0}$ versus $H_1:\textbf{R}\boldsymbol{\beta}\not \le \textbf{0}$,
where ``$\not\le$'' implies that at least one inequality goes the other direction. Third, we briefly consider the three hypothesis test $H_0:\textbf{R}\boldsymbol{\beta}=\textbf{0}$ versus $H_1:\textbf{R}\boldsymbol{\beta}\le\textbf{0}$ (with $\textbf{R}\boldsymbol{\beta}=\textbf{0}$ excluded) versus $H_2:\textbf{R}\boldsymbol{\beta}\not \le\textbf{0}$.

The model is reparametrized so that the linear combination of parameters of interest, i.e., $\vet{\theta}=\textbf{R}\bm\beta$, is perpendicular to the nuisance parameters, i.e., $\boldsymbol{\gamma}=\textbf{D}\bm\beta$, i.e.,
\[
\left[\begin{array}{c}\vet{\theta}\\ \vet{\gamma}\end{array}\right]
=
\left[\begin{array}{c}\textbf{R}\\ \textbf{D}\end{array}\right]\vet{\beta}
=\textbf{T}\boldsymbol{\beta},
\]
where the $r_2\times K$ matrix $\textbf{D}$ contains $r_2=K-r_1$ independent rows of $\textbf{P}_{\textbf{R}}^{\perp}\textbf{X}_{\boldsymbol{\beta}}'\vet{\Sigma}^{-1}\textbf{X}_{\boldsymbol{\beta}}$, where the orthogonal projection matrix is given by $\textbf{P}_{\textbf{R}}^{\perp}=\textbf{I}_K-\textbf{R}'\left(\textbf{R}\textbf{R}'\right)^{-1}\textbf{R}$. Subsequently, the model can be written as
\[
\textbf{y}=
\textbf{X}_{\boldsymbol{\theta}}\boldsymbol{\theta}+\textbf{X}_{\boldsymbol{\gamma}}\vet{\gamma}+\vet{\epsilon},
\]
where $\textbf{X}_{\boldsymbol{\theta}}$ contain the first $r_1$ columns of $\textbf{X}\textbf{T}^{-1}$ that are regressed on $\boldsymbol{\theta}$ and $\textbf{X}_{\boldsymbol{\gamma}}$ contains the remaining $r_2$ columns of $\textbf{X}\textbf{T}^{-1}$ that are regressed on $\boldsymbol{\gamma}$. The null hypothesis can then be written as $H_0:\vet{\theta}=\textbf{0}$ versus $H_1:\vet{\theta}\in\mathbb{R}^{r_1}$ and the nonnested hypothesis test can be written as $H_0:\vet{\theta}\le\textbf{0}$ versus $H_1:\vet{\theta}\not \le \textbf{0}$. Further note that the ML estimates of $\boldsymbol{\theta}$ and $\boldsymbol{\gamma}$ are independent because $\left([\textbf{X}\textbf{T}^{-1}]'\boldsymbol{\Sigma}^{-1}[\textbf{X}\textbf{T}^{-1}]\right)^{-1}=\mbox{diag}\left(\left(\textbf{X}_{\boldsymbol{\theta}}'\boldsymbol{\Sigma}^{-1}\textbf{X}_{\boldsymbol{\theta}}\right)^{-1},\left(\textbf{X}_{\boldsymbol{\gamma}}'\boldsymbol{\Sigma}^{-1}\textbf{X}_{\boldsymbol{\gamma}}\right)^{-1}\right)$ which is a direct consequence of the choice of $\textbf{D}$.

Throughout this paper, the free parameters under a hypothesis have a hypothesis index to make it explicit that the parameters under different hypotheses have different interpretations and therefore different priors. For example, the population variances under $H_0$ and $H_1$ are denoted by $\sigma_0^2$ and $\sigma_1^2$, respectively.

\section{Testing a Precise Hypothesis}\label{nullHypo}
The following definition will be used for information inconsistency when testing a precise hypothesis.

\begin{definition}
\label{incon1}
A Bayes factor, $B_{10}$, is called information inconsistent for testing $H_0:\vet{\theta}=\textbf{0}$ versus $H_1:\vet{\theta}\not=\textbf{0}$ if
there exists a sequence of possible data,  with corresponding $\hat{\vet{\theta}}_i$ that satisfy $|\hat{\vet{\theta}}_i|\rightarrow\infty$ as $i \rightarrow \infty$, for which the Bayes factors $B_{10} < B_{10}^*<\infty$.
\end{definition}

\subsection{Conjugate priors}
In the conjugate case, the conditional prior of $\bm \theta \mid \sigma^2_1$ under $H_1$ has a multivariate normal distribution and the marginal prior of $\sigma^2_t$, for $t=0$ or 1, has a scaled inverse chi-squared distribution, resulting in
\begin{eqnarray}
\nonumber \pi_1(\vet{\theta},\vet{\gamma}_1,\sigma^2_1)&= &\pi_1(\vet{\theta}  \mid \sigma^2_1)\times \pi_1(\bm\gamma_1)\times \pi_1(\sigma^2_1)\\
\label{conjPrior1} &\propto&N_{\boldsymbol{\theta} \mid \sigma^2_1}(\textbf{0},\sigma_1^2 \bm \Omega)\times 1\times \mbox{inv-}\chi^2_{\sigma^2_1}(s_1^2,\nu_1)\\
\nonumber \pi_0(\vet{\gamma}_0,\sigma_0^2) &=& \pi_0(\vet{\gamma}_0)\times\pi_0(\sigma_0^2)\\
\label{conjPrior0}&\propto& 1\times\mbox{inv-}\chi^2_{\sigma^2_0}(s_0^2,\nu_0) \,.
\end{eqnarray}
The scaled inverse chi-squared distribution is used (instead of the inverse gamma distribution) because of the natural relation between the prior degrees of freedom $\nu_t$ and the sample size $n$ \citep{Gelman}. When setting the prior degrees of freedom equal to $\nu_t=0$, we obtain the objective improper prior, $\pi_t(\sigma_t^2)\propto\sigma_t^{-2}$, for $t=0$ or 1, and when additionally setting $\bm \Omega=g\left(\textbf{X}_{\bm\theta}'\bm \Sigma^{-1}\textbf{X}_{\bm\theta}\right)^{-1}$, we obtain Zellner's $g$-prior.

Denoting the ML estimate of $\bm{\theta}$ by $\hat{\boldsymbol{\theta}}= \left(\textbf{X}_{\bm\theta}'\boldsymbol{\Sigma}^{-1}\textbf{X}_{\bm\theta}\right)^{-1}\textbf{X}_{\bm\theta}'\boldsymbol{\Sigma}^{-1}\textbf{y}$ and the sums of squares by $s^2_{\textbf{y}}=(\textbf{y}-\textbf{X}_{\bm\theta}\hat{\boldsymbol{\theta}}-\textbf{X}_{\bm\gamma}\hat{\boldsymbol{\gamma}})'\bm\Sigma^{-1}(\textbf{y}-\textbf{X}_{\bm\theta}\hat{\boldsymbol{\theta}}-\textbf{X}_{\bm\gamma}\hat{\boldsymbol{\gamma}})$, a standard calculation yields that the Bayes factor of $H_1$ against $H_0$, based on the conjugate priors in \eqref{conjPrior1} and \eqref{conjPrior0}, is
\begin{eqnarray}
B_{10} &=& C_1 \times \frac{\left(s_1^2\nu_1+s^2_{\textbf{y}}+\hat{\boldsymbol{\theta}}'
\left(\left(\textbf{X}_{\bm\theta}'\boldsymbol{\Sigma}^{-1}\textbf{X}_{\bm\theta}\right)^{-1}+\bm \Omega \right)^{-1}
\hat{\boldsymbol{\theta}}\right)^{-(n+\nu_1-r_2)/2}}
{\left(s_0^2\nu_0+s^2_{\textbf{y}}+\hat{\boldsymbol{\theta}}'
\textbf{X}_{\bm\theta}'\boldsymbol{\Sigma}^{-1}\textbf{X}_{\bm\theta}\hat{\boldsymbol{\theta}}\right)^{-(n+\nu_0-r_2)/2}}\,,
\label{B10nullConj}
\end{eqnarray}
where the constant is
$$C_1 = \frac{(\nu_1/2)^{\nu_1/2}s_1^{\nu_1}\Gamma\left(\frac{\nu_0}{2}\right)\Gamma\left(\frac{n+\nu_1-r_2}{2}\right)}{(\nu_0/2)^{\nu_0/2}s_0^{\nu_0}\Gamma\left(\frac{\nu_1}{2}\right)\Gamma\left(\frac{n+\nu_0-r_2}{2}\right)} 2^{(\nu_1-\nu_0)/2}
|{\bm \Omega} +(\textbf{X}_{\bm\theta}'\boldsymbol{\Sigma}^{-1}\textbf{X}_{\bm\theta})^{-1}|^{-\frac{1}{2}}
|\textbf{X}_{\bm\theta}'\boldsymbol{\Sigma}^{-1}\textbf{X}_{\bm\theta}|^{-\frac{1}{2}} \,.
$$
The following result is immediate.
\begin{lemma}\label{lemma1}
As $|\hat{\boldsymbol{\theta}}| \rightarrow \infty$, the Bayes factor in \eqref{B10nullConj} satisfies
$ B_{10} \rightarrow   0$ if $\nu_0 < \nu_1$; $ B_{10} \rightarrow  \infty $ if $\nu_0 > \nu_1$; and, if $\nu_0 = \nu_1$,
$$
B_{10} \leq   C_1  \left(  \limsup_{ |\hat{\boldsymbol{\theta}}| \rightarrow \infty} \frac{\hat{\boldsymbol{\theta}}'
\textbf{X}_{\bm\theta}'\boldsymbol{\Sigma}^{-1}\textbf{X}_{\bm\theta}\hat{\boldsymbol{\theta}}}
{\hat{\boldsymbol{\theta}}'
\left(\left(\textbf{X}_{\bm\theta}'\boldsymbol{\Sigma}^{-1}\textbf{X}_{\bm\theta}\right)^{-1}+\bm \Omega \right)^{-1}
\hat{\boldsymbol{\theta}}}\right)^{\frac{(n+\nu-r_2)}{2}}  = \ C_1 \ \left(1+\lambda_{max}\right)^{(n+\nu-r_2)/2} < \infty \,,
$$
where $\lambda_{max}$ is the largest eigenvalue of $\textbf{X}_{\bm\theta}'\boldsymbol{\Sigma}^{-1}\textbf{X}_{\bm\theta}\bm \Omega$.
\end{lemma}

\begin{remark}
Setting $\nu_0 < \nu_1$ seems logical because it implies that the prior for $\sigma^2_1$ is more concentrated than the prior for $\sigma_0^2$ (consistent with a nonzero mean explaining some of the variation compared to a zero mean). This choice however results in a disastrously information inconsistent Bayes factor, with the conclusion being that the null hypothesis is certainly true when $|\hat{\boldsymbol{\theta}}| \rightarrow \infty$.
\end{remark}

\begin{remark}
Setting $\nu_0 = \nu_1$ is the usual choice, which still results in an information inconsistent Bayes factor. Note that the prior degrees of prior would be set to 0 in the objective Bayesian approach. The impact of this inconsistency will be discussed below for the special case of the univariate $t$-test.
\end{remark}

\begin{remark}
Setting $\nu_0 > \nu_1$ would not be a logical choice because the prior for $\sigma_0^2$ is more concentrated than the prior for $\sigma_1^2$, even though the regression coefficients $\bm\theta$ under $H_0$ are restricted to 0. The resulting Bayes factor, however, is information consistent. A special case of this choice arises from setting the prior for the variance under $H_0$ to be proportional to the conditional prior of the variance given $\bm\theta=\textbf{0}$ under $H_1$, i.e., $\pi_0(\sigma^2)=\pi_1(\sigma^2 \mid \vet{\theta}=\textbf{0})=\mbox{inv-}\chi^2(\tfrac{\nu_1}{\nu_1+r_1}s_1^2,\nu_1+r_1)$, so that $\nu_0=\nu_1+r_1$. The Bayes factor can then be expressed as the Savage-Dickey density ratio \citep{Dickey:1971}, $B_{10}=\frac{\pi_1(\bm\theta=\textbf{0})}{\pi_1(\bm\theta=\textbf{0}|\textbf{X})}$, where the marginal prior and the posterior of $\bm\theta$ have a Student $t$-distribution.
\end{remark}

\subsubsection{Practical implications for a univariate test under dependence}
The practical importance of information inconsistency is explored for the objective prior with $\nu_1=\nu_0=0$ for a univariate t-test of $H_0:\theta=0$ versus $H_1:\theta\not =0$ with correlated data. Specifically, consider $r_1=1$, $r_2=0$, $\textbf{X}_{\bm\theta}=\textbf{1}_n$, and $\bm \Omega=1$, with $\bm \Sigma$ being the correlation matrix with identical correlations $\rho$ in the off-diagonal elements. The t-statistic, $t= \frac{\hat{\theta}  \sqrt{ \textbf{1}_n' \bm\Sigma^{-1}\textbf{1}_n}} {s_{\textbf{y}}/\sqrt{n-1}}$, then has a $t$-distribution with $n-1$ degrees of freedom under $H_0$. The Bayes factor in \eqref{B10nullConj} can then be expressed as a function of the t-statistic, namely
\[
B_{10} = \left(1+\frac{n}{1+(n-1)\rho}\right)^{-1/2} \left(1 - \frac{nt^2}{[t^2 +n-1][n+1+(n-1)\rho]}\right)^{-n/2}.
\]
The limiting value of the Bayes factor, as $|t|$ goes to infinity, is
\begin{eqnarray*} \lim_{|t| \rightarrow \infty} B_{10} &=& \left(1+\frac{n}{1+(n-1)\rho}\right)^{-1/2} \left(1 - \frac{n}{n+1+(n-1)\rho}\right)^{-n/2} \\
&=& \left\{
      \begin{array}{ll}
        (1+n)^{(n-1)/2}, & \hbox{if \ $\rho=0$;} \\
        \left(1+\frac{2n}{n+1}\right)^{-1/2} \left(\frac{3n+1}{n+1}\right)^{n/2} \approx 3^{(n-1)/2}, & \hbox{if \ $\rho=0.5$;} \\
        2^{(n-1)/2}, & \hbox{if \ $\rho=1$.}
      \end{array}
    \right.
\end{eqnarray*}
Hence, the correlation can dramatically affect the situation. Table \ref{table1} provides the limiting value of the Bayes factor as $|t|$ goes to $\infty$ for different choices of the correlation $\rho$ and different sample sizes varying from $n=2$ to a sample size of $n=20$. The table also provides the Bayes factor when $t=4$ to check whether inconsistency starts coming into play for a large $t$-value. As comparisons, the corresponding two-sided $p$-values are also provided, as well as the upper bound $B_{10} < 1/[-ep \log p]$, which is a bound over a large nonparametric class of priors (derived in \cite{Sellke:2001}).

When there is zero correlation, the limit $(n+1)^{(n-1)/2}$ is large for sample sizes larger than 6, so that information inconsistency is then not problematical practically. For large correlations on the other hand, and especially when $\rho$ is close to 1, the limiting values can be quite small, arguing against the use of objective conjugate priors.

\begin{figure}[t]
\centering
\makebox{\includegraphics[width=15.5cm,height=4.7cm]{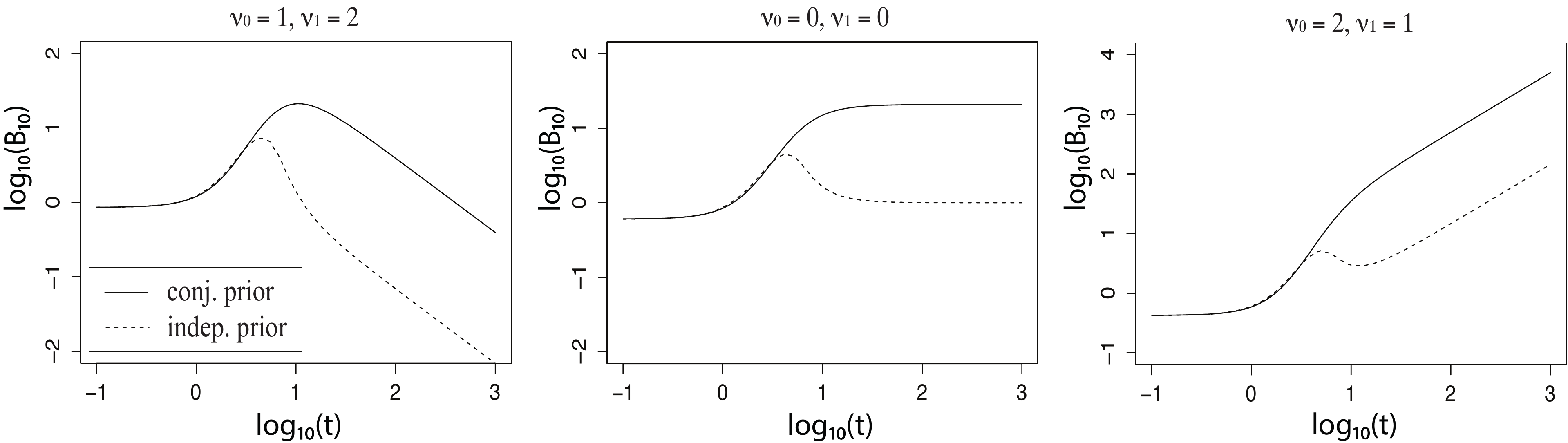}}
\caption{The Bayes factor $B_{10}$ based on the conjugate prior (solid line) and independence prior (dashed line) as a function of $t$-values when $n=7$, $\rho=.5$, $s_{\textbf{y}}^2=n-1=6$, $s_0^2=s_1^2=1$, and different choices for the prior degrees of freedom $\nu_0$ and $\nu_1$.} \label{fig1}
\end{figure}

Figure \ref{fig1} displays the Bayes factor as a function of $\log_{10}(t)$ when using conjugate priors (solid lines) and $n=7$, $\rho=.5$, $s_{\textbf{y}}^2=n-1=6$, $s_0^2=s_1^2=1$, and different choices for the prior degrees of freedom, namely $(\nu_0,\nu_1)=(0,0)$, $(1,2)$ or $(2,1)$. As can be seen, if $\nu_0=\nu_1=0$, the logarithm of the Bayes factor converges to $\log_{10}(20.8)=1.32$ (Table \ref{table1}); if $\nu_0>\nu_1$. Furthermore, if $\nu_0<\nu_1$ (or $\nu_0>\nu_1$), the evidence goes to $\infty$ for $H_0$ (or $H_1$) as $t\rightarrow\infty$ implying information inconsistency (or information consistency).

\begin{table}[t]
\caption{Limiting values of the Bayes factor for a univariate $t$-test as $|t|\rightarrow\infty$ for different choices of the sample size $n$ and the correlation $\rho$. Additionally Bayes factors and two-sided p-values are given when $t=4$. The approximation $1/[-ep \log p]$ is an upper bound of the evidence against $H_0$ \citep{Sellke:2001}.}
{\small \begin{tabular}{llccccc}
  \hline
& $n$ & 2 & 5 &7 & 10 &20 \\
\hline
\hline
$\rho=0$ & \mbox{limit}  & 1.73 & 36 & 512 & $4.85 \times 10^4$ & $1.79 \times 10^{11}$ \\
& $B_{10}$ for $t=4$ & 1.55 & 6.36 & 12.21 & 23.61 & 66.20 \\
  \hline
$\rho=0.5$ & \mbox{limit}  & 1.53 & 7.10 & 20.8 & 106 & $2.01 \times 10^4$ \\
& $B_{10}$ for $t=4$ & 1.42 & 3.46 & 5.31 & 8.54 & 20.71 \\
  \hline
$\rho \approx 1$ & \mbox{limit}  & 1.41 & 4 & 8 & 22.6 & 724 \\
& $B_{10}$ for $t=4$ & 1.34 & 2.76 & 3.44 & 4.86 & 9.47 \\
  \hline
\hline
& $p$-value for $t=4$ & 0.156 & 0.016 & 0.0071 & 0.0031 & 0.00077 \\
  \hline
& $1/[-ep \log p]$  & 2.25 & 7.81 & 13.47 & 24.40 & 72.01 \\
  \hline
\end{tabular}}
\label{table1}
\end{table}

\subsection{Mixtures of conjugate priors}

Although use of conjugate priors in testing is common, it has long been argued (starting with \cite{Jeffreys}) that fatter tailed prior distributions should be used. One such class that is increasingly popular is the class of scale mixtures of conjugate priors. This class results in information consistent Bayes factors if the prior on $g$ is thick enough, as shown by the following lemmas which generalize the result in \cite{Liang:2008} for $\nu_0 = \nu_1 = 0,$ $\bm{\Sigma} = \bm{I},$ and $\bm{\Omega} = g (\textbf{X}_{\boldsymbol{\theta}}' \bm\Sigma^{-1} \textbf{X}_{\boldsymbol{\theta}})^{-1}$.

\begin{lemma} \label{genmixt}
 Let $\bm{\theta} \mid g, \bm{\gamma}_1,  \sigma_1^2 \sim N_{r_1}(\bm{0}_{r_1}, g \sigma_1^2 \Om )$, where $\sigma_1^2$ has the prior specified in  (3) and $g$ has a prior with density $\pi(g)$. If $\nu_0 > \nu_1$, any $\pi(g)$ with positive support yields an information-consistent $B_{10}$. The condition
$$
\int_0^\infty (g+1)^{(n-r_1-r_2+\nu_1)/2} \pi(g) \mathrm{d} g= \infty
$$
is necessary and sufficient for information consistency whenever $\nu_0 = \nu_1$, and necessary whenever $\nu_0 < \nu_1$.
\end{lemma}

\noindent \textit{Proof:} See Appendix~\ref{appA}.

The maximum number of finite moments that the prior on $g$ can have to achieve information consistency increases with the sample size $n$ and decreases with the number of predictors $K = r_1+r_2$. Lemma~\ref{genmixt} gives us a complete description for all scale mixtures of conjugate priors whenever $\nu_0 \ge \nu_1,$ but only gives us a necessary condition for information consistency for $\nu_0 < \nu_1$. The lemma below characterizes the behavior of polynomial-tailed priors on $g$ in this latter case, and provides partial results for priors with thinner- and thicker-than-polynomial priors on $g$.

\begin{lemma} \label{tail}
  Suppose $\nu_0 < \nu_1$ and let $\bm{\theta} \mid g, \bm{\gamma}_1, \sigma^2_1 \sim N_{r_1}(\bm{0}_{r_1}, g \sigma^2_1 \Om )$, where $\sigma^2_1$ has the prior specified in (3) and $g$ has a prior with density $\pi(g)$. Then, the following are true:
  \begin{enumerate}
      \item If there exist $0 < M < \infty$ and $0 < K < \infty$ such that for all $g \ge M$,  $\pi(g) \ge K g^{-\alpha}$ for $\alpha >1$, $B_{10}$ is information consistent whenever $\alpha < (n-r_1-r_2+\nu_0)/2+1$.
      \item If there exist $0 < M' < \infty$ and $0 < K' < \infty$ such that for all $g \ge M'$,  $\pi(g) \le K' g^{-\alpha}$ for $\alpha >1$, $B_{10}$ is information inconsistent whenever $\alpha \ge  (n-r_1-r_2+\nu_0)/2+1$.
  \end{enumerate}
  [NB: All of the priors on $g$ considered in \cite{Liang:2008} satisfy both conditions.]
\end{lemma}

\noindent \textit{Proof:} See Appendix~\ref{appB}.

Note that the Zellner-Siow prior \citep{ZellnerSiow:1980} (which was the first proposed information consistent prior for this situation) and the Hyper-$g$ prior \citep{Liang:2008} satisfy both conditions because they have polynomial tails.

\subsection{Independence priors}

\subsubsection{Semi-conjugate prior}
A feature of the conjugate prior that is sometimes questioned is the dependence induced between $\vet{\theta}$ and $\sigma^2$; in objective Bayesian analysis this is hard to avoid (only $\sigma$ is available to provide an objective scale for $\vet{\theta}$), but it does seem rather
arbitrary. Hence it is of interest to also investigate information consistency using semi-conjugate priors of the form
\begin{eqnarray*}
\pi_1(\vet{\theta},\vet{\gamma}_1,\sigma_1^2)&=&\pi_1(\vet{\theta})\times\pi_1(\vet{\gamma}_1)\times\pi_1(\sigma_1^2)\\
&\propto& N_{\boldsymbol{\theta}}(\textbf{0},\bm \Omega)\times 1 \times\mbox{inv-}\chi^2_{\sigma^2_1}(s_1^2,\nu_1)\\
\pi_0(\vet{\gamma}_0,\sigma_0^2)&=&\pi_0(\vet{\gamma}_0)\times\pi_0(\sigma_0^2)\\
&\propto& 1 \times\mbox{inv-}\chi^2_{\sigma^2_0}(s_0^2,\nu_0).
\end{eqnarray*}
With these semi-conjugate priors, the Bayes factor becomes
\begin{eqnarray}
B_{10} &=& C_2 \times\frac{\int \left(\nu_1 s_1^2+s^2_{\textbf{y}}+(\boldsymbol{\theta}-\hat{\boldsymbol{\theta}})'\textbf{X}'_{\bm\theta}\boldsymbol{\Sigma}^{-1}\textbf{X}_{\bm\theta}(\boldsymbol{\theta}-\hat{\boldsymbol{\theta}}) \right)^{-\frac{n-r_2+\nu_1}{2}}N_{\boldsymbol{\theta}}(\textbf{0},\bm\Omega) d\boldsymbol{\theta}}{\left(\nu_0 s_0^2+s^2_{\textbf{y}}+\hat{\boldsymbol{\theta}}'\textbf{X}'_{\bm\theta}\boldsymbol{\Sigma}^{-1}\textbf{X}_{\bm\theta}\hat{\boldsymbol{\theta}} \right)^{-\frac{n-r_2+\nu_0}{2}}} \,,
\label{BF10nullIndep}
\end{eqnarray}
where
$$ C_2 =\frac{(\nu_1/2)^{\nu_1/2}s_1^{\nu_1}\Gamma\left(\frac{\nu_0}{2}\right)\Gamma\left(\frac{n+\nu_1-r_2}{2}\right)}{(\nu_0/2)^{\nu_0/2}s_0^{\nu_0}\Gamma\left(\frac{\nu_1}{2}\right)\Gamma\left(\frac{n+\nu_0-r_2}{2}\right)} 2^{(\nu_1-\nu_0)/2} \,.
$$

\begin{lemma}
\label{lemma.indconj}
As $|\hat{\boldsymbol{\theta}}| \rightarrow \infty$, the Bayes factor in \eqref{BF10nullIndep}, based on the independent semi-conjugate prior, behaves as follows:
\begin{eqnarray*}
B_{10} &\rightarrow& \left\{
                \begin{array}{ll}
                  0 & \hbox{if $\nu_0 < \nu_1$;} \\
         1 & \hbox{if $\nu_0 = \nu_1$;} \\
                  \infty & \hbox{if $\nu_0 > \nu_1$.}
                \end{array}
              \right.
\end{eqnarray*}
\end{lemma}
\textit{Proof: See Appendix \ref{proof.indconj}.}

Note that, in the typical case of $\nu_0=\nu_1$, we observe an even worse case of information inconsistency than for the conjugate prior because the relative evidence between $H_1$ and $H_0$ goes to 1 when there appears to be overwhelming evidence for $H_1$; in contrast, for the conjugate prior case, the limiting Bayes factor -- while nonzero -- was at least exponentially small in $n$.

The intuition behind this result is that very large $\hat{\bm\theta}$ are equally unlikely under $H_1$ and $H_0$, due to the light-tailed normal prior for $\bm\theta$ under $H_1$. Furthermore, the limits are the same as in the conjugate case if $\nu_0\not =\nu_1$. Hence, the choice of the prior degrees of freedom plays a crucial role in information inconsistency, even when the variance is apriori independent of $\bm{\theta}$.

Figure \ref{fig1} also displays the Bayes factor, based on the independence prior, as a function of $\log_{10}(t)$ for the univariate $t$-test when the data correlation is $\rho=.5$ (dashed line). As can be seen the Bayes factor based on the independence prior and the conjugate prior with the same hyper parameters are approximately equal for absolute $t$ values smaller than approximately $\log_{10}(.5)$. For larger t values, the flatter tails of the independence priors start to have an effect resulting in a decrease of the Bayes factor, relative to the Bayes factor based on the conjugate priors.

\subsubsection{Fatter-tailed independence priors}
It is somewhat unfair to use an independent normal prior for model comparison here since, from \cite{Jeffreys}, the use of fatter tailed priors has been recommended. To keep the discussion of fatter tailed priors simple, we consider only the one dimensional case (i.e., $r_2=0$), and restrict the prior $\pi_1(\theta)$ to be
a $t$-distribution with mean 0, scale $\tau$ (fixed) and degrees of freedom $\nu$, i.e.,
$$\pi_1(\theta)=\frac{\Gamma((\nu+1)/2)}{\sqrt{\nu\pi} \ \Gamma(\nu/2)\tau}\left(1+\frac{\theta^2}{\nu \tau^2}\right)^{-\frac{\nu+1}{2}} \,.$$
Then Theorem 3.3 in \cite{fan1992behaviour} shows that,
as $|\hat \theta| \rightarrow \infty$,
$$
B_{10} = C \frac{\frac{\Gamma((n^*+1)/2)}{\sqrt{n^*\pi} \ \Gamma(n^*/2)\sqrt{V}} \left(1+\frac{\hat{{\theta}}^2}
{n^* V} \right)^{-\frac{n^*+1}{2}}+\frac{\Gamma((\nu+1)/2)}{\sqrt{\nu\pi} \ \Gamma(\nu/2)\tau}\left(1+\frac{\hat{\theta}^2}{\nu \tau^2}\right)^{-\frac{\nu+1}{2}} }{\left(\nu_0 s_0^2+s^2_{\textbf{y}}+\hat{{\theta}}'\textbf{X}'_{\theta}\boldsymbol{\Sigma}^{-1}\textbf{X}_{\theta}\hat{{\theta}} \right)^{-\frac{n+\nu_0}{2}}} \times (1+o(1))\,,
$$
where $n^* = n+\nu_1-1$, $ V = (\nu_1 s_1^2+s^2_{\textbf{y}})/[n^*\textbf{X}'_{\theta}\boldsymbol{\Sigma}^{-1}\textbf{X}_{\theta}]$ and
$$C = \frac{(\nu_1/2)^{\nu_1/2}s_1^{\nu_1}\sqrt{n^*\pi} \ \Gamma(n^*/2)\sqrt{V}}{\Gamma\left(\nu_1/2\right) \ (\nu_1 s_1^2+s^2_{\textbf{y}})^{(n+\nu_1) / 2}} \,.$$
Thus, as $|\hat \theta| \rightarrow \infty$,
\begin{eqnarray*}
B_{10} &\rightarrow& \left\{
                \begin{array}{ll}
                  0 & \hbox{if $n+\nu_0 < \min\{n-1+\nu_1,\nu+1\}$;} \\
         \hbox{constant} & \hbox{if $n+\nu_0 = \min\{n-1+\nu_1,\nu+1\}$;} \\
                  \infty & \hbox{if $n+\nu_0 > \min\{n-1+\nu_1,\nu+1\}$.}
                \end{array}
              \right.
\end{eqnarray*}
Since $n \geq 2$, if $0 < \nu < 1$ it will be true that $n+\nu_0 > \min\{n-1+\nu_1,\nu+1\}$ so that $B_{10}$ will be information consistent. For the commonly used Cauchy prior ($\nu=1$), information consistency also holds, except for the case when $n=2$ and $\nu_0=0$ (this last corresponding to the objective prior for $\sigma_1^2$).
It is interesting that information consistency does hold for this last case when $\pi_1(\theta)$ is chosen to be $\mbox{Cauchy}(0, \sigma_1)$ \citep[cf.][]{Liang:2008} and $\nu_1 =0$; thus, once again, insisting on prior independence
of $\sigma_1^2$ and $\theta$ only appears to worsen the problem of information inconsistency.

\subsection{Adaptive priors}
Another approach to Bayesian hypothesis testing is to let the prior under $H_1$ adapt to the likelihood, as in \cite{george2000calibration} and \cite{hansen2001model}.

\begin{example}
For the $g$-prior in the $t$-test, when the $t$-statistic $t=\sqrt{\frac{\hat{\bm\theta}'\textbf{X}_{\bm\theta}'\bm\Sigma^{-1}\textbf{X}_{\bm\theta}\hat{\bm\theta}}{s_{\textbf{y}}^2/(n-1)}} >1$, the marginal likelihood under $H_1$ is maximized (calculus) for the choice $g=\frac{n-r_2-r_1}{r_1(n-1)}t^2-1$. The Bayes factor for this choice equals
\begin{eqnarray*}
B_{10}&=&\left(\frac{r_1(n-1)}{t^2(n-r_1-r_2)}\right)^{\frac{r_1}{2}}\left(\frac{(n-1+t^2)(n-r_1-r_2)}{(n-1)(n-r_2)}\right)^{\frac{n-r_2}{2}},
\end{eqnarray*}
which is information consistent. For a univariate t test, with $r_1=1$ and $r_2=0$, the resulting Bayes factor can be expressed as $B_{10}=\frac{1}{|t|}\left(\frac{n-1+t^2}{n}\right)^{\frac{n}{2}}$.
\end{example}

The following lemma generalizes the result in \cite{Liang:2008} for $\nu_0 = \nu_1 = 0,$ $\bm{\Sigma} = \bm{I},$ and $\bm{\Omega} = g (\textbf{X}_{\boldsymbol{\theta}}' \bm\Sigma^{-1} \textbf{X}_{\boldsymbol{\theta}})^{-1}$.

\begin{lemma} \label{adaptive}
 Let $\bm{\theta} \mid g, \bm{\gamma}_1, \sigma^2_1 \sim N_{r_1}(\bm{0}_{r_1}, g \sigma^2_1 \Om )$, where $\sigma^2_1$ has the prior specified in (3). If $g > 0$ is set by maximizing $B_{10}$, information  consistency holds.
\end{lemma}

\noindent \textit{Proof:} See Appendix~\ref{appD}.

Lemma~\ref{adaptive} establishes information consistency for all $\nu_0$ and $\nu_1$. This is in contrast with the results in previous sections, where the behavior of $B_{10}$ depends (sometimes rather strongly) on $\nu_0$ and $\nu_1$.



\section{One-Sided Hypothesis Testing}\label{nonnestedHypo}
The following definition will be used for information consistency for a one-sided testing problem.

\begin{definition}
\label{incon2}
A Bayes factor is {\em information consistent}, for a one-sided hypothesis test of $H_0:\bm\theta\le\textbf{0}$ versus $H_1:\bm\theta\not\le\textbf{0}$, if $B_{10}\rightarrow\infty$ as $|\hat{\bm\theta}| \rightarrow \infty$ with at least one coordinate of $\hat{\bm\theta}$ going to $\infty$, and $B_{10}\rightarrow 0$, as all coordinates of $\hat{\bm\theta}$ go to $-\infty$. If this does not hold, the Bayes factor is called information inconsistent.
\end{definition}

We shall denote the subspaces under $H_0$ and $H_1$ as $\bm\Theta_0=\{\bm\theta \mid \bm\theta\le\textbf{0}\}$ and $\bm\Theta_1=\{\bm\theta \mid \bm\theta\not\le\textbf{0}\}$, respectively.

\subsection{Conjugate prior}
When testing nonnested hypotheses, it is common to formulate an encompassing prior $\pi$ on the joint space $\bm{\Theta}=\bm{\Theta}_0\cup\bm{\Theta}_1$ and specify truncations of this prior under $H_0$ and $H_1$ \citep[e.g.,][]{Berger:1999,Klugkist:2007}. As in the null hypothesis test, the encompassing conjugate prior is centered on the boundary of the subspaces under investigation, i.e.,
\begin{eqnarray}
\pi(\bm{\theta},\bm{\gamma},\sigma^2)\propto N_{\boldsymbol{\theta} \mid \sigma^2}(\textbf{0},\sigma^2\bm\Omega)\times \mbox{inv-}\chi^2_{\sigma^2}(s^2,\nu),
\label{conjprior2}
\end{eqnarray}
with a flat improper prior for $\bm{\gamma}$.
The priors under the nonnested hypotheses $H_t$, for $t=0$ or 1, can then be expressed as
\begin{eqnarray}
\pi_t(\bm{\theta} \mid \sigma^2) &=& \pi(\bm{\theta} \mid \sigma^2)I_{\boldsymbol{\Theta}_t}(\boldsymbol{\theta})/P_{\pi}(\bm{\theta}\in\bm{\Theta}_t \mid \sigma^2) \,,
\label{priornonnested}
\end{eqnarray}
$\pi_t(\sigma^2)=\pi(\sigma^2)$, and $\pi_t(\bm{\gamma})=\pi(\bm{\gamma})$, with the denominator in \eqref{priornonnested} being equal to the conditional prior probability of $\bm\Theta_t$ under the joint prior on $\bm\Theta$, i.e., $P_{\pi}(\bm{\theta}\in\bm{\Theta}_t \mid \sigma^2)=\int_{\bm\Theta_t}N_{\boldsymbol{\theta} \mid \sigma^2}(\textbf{0},\sigma^2\bm\Omega)d\bm\theta>0$.

The Bayes factor for the one-sided hypothesis test based on the conjugate priors can then be expressed as
\begin{equation}
B_{10} = \left(P_{\pi}(\bm{\theta}\le\textbf{0} \mid \sigma^2=1)^{-1}-1\right)^{-1} \left(P_{\pi}(\bm{\theta}\le\textbf{0} \mid \textbf{y})^{-1}-1\right).
\label{B10onesided}
\end{equation}
The derivation is similar to that in \cite{Mulder:2014a}. 
The prior and posterior probabilities that the constraints hold under the encompassing model can be computed as the proportion of draws satisfying the constraints. Also note that the conditional prior probability of $\bm{\theta}\le\textbf{0}$ is completely determined by the prior covariance matrix $\bm\Omega$ and is independent of $\sigma^2$ (therefore we can set $\sigma^2=1$ in \eqref{B10onesided}). This is a direct result of centering the encompassing prior on the point of interest $\textbf{0}$. For example, if $\bm\Omega=\textbf{I}_{r_1}$, then $P_{\pi}(\vet{\theta}\le\textbf{0} \mid \sigma^2)=2^{-r_1}$, $\forall \sigma^2>0$. In the $g$-prior with $\bm\Omega=g\sigma^{2}(\textbf{X}_{\bm\theta}\bm\Sigma^{-1}\textbf{X}_{\bm\theta})^{-1}$, the prior probability is completely determined by the covariance structure of the predictors.

As can be concluded from \eqref{B10onesided}, a Bayes factor for a one-sided hypothesis test is information consistent if $P_{\pi}(\bm{\theta}\le\textbf{0} \mid \textbf{y})\rightarrow 0$ as $|\hat{\bm\theta}| \rightarrow \infty$ with at least one coordinate of $\hat{\bm\theta}$ going to $\infty$, and $P_{\pi}(\bm{\theta}\le\textbf{0} \mid \textbf{y})\rightarrow 1$ as all coordinates of $\hat{\bm\theta}$ go to $-\infty$.

\begin{lemma}
\label{inconBF2}
$P_{\pi}(\bm{\theta}\le\textbf{0} \mid \textbf{y})$ is bounded away from 0 and 1 for all $\textbf{y}$. Hence
$B_{10}$ is information inconsistent.

If $\hat{\bm\theta} = c \bm v$ and $c \rightarrow \infty$, then
$$P_{\pi}(\bm{\theta}\le\textbf{0} \mid \textbf{y}) \rightarrow P_{\pi}(\bm{\xi}\le\textbf{0} \mid \textbf{y}) \,,$$
where $\bm{\xi}$ has a multivariate t distribution with mean
$${\bm v}^* = \frac{(\textbf{X}_{\bm\theta}'\bm{\Sigma}^{-1}\textbf{X}_{\bm\theta}+
\bm\Omega^{-1})^{-1}\textbf{X}_{\bm\theta}'\bm{\Sigma}^{-1}\textbf{X}_{\bm\theta}{\bm v}}
{(n+\nu-r_2)^{-1/2}({\bm v}'((\textbf{X}_{\bm\theta}'\bm{\Sigma}^{-1}\textbf{X}_{\bm\theta})^{-1}+
\bm\Omega)^{-1}{\bm v})^{1/2}}\,,$$
 scale matrix $(\textbf{X}_{\bm\theta}'\bm{\Sigma}^{-1}\textbf{X}_{\bm\theta}+\bm\Omega^{-1})^{-1}$, and  $n+\nu-r_2$ degrees of freedom.

\end{lemma}
\textit{Proof: See Appendix \ref{proof.inconBF2}.}

\subsubsection{Practical implications for a univariate one-sided test under dependence}
We investigate the practical importance of information inconsistency for a univariate one-sided $t$-test under dependence of $H_0:\theta\le 0$ versus $H_1:\theta>0$, with $\nu=0$, $r_1=1$, $r_2=0$, $\textbf{X}_{\bm\theta}=\textbf{1}$, $\bm\Omega=1$, and $\bm\Sigma=\rho\textbf{J}_n+(1-\rho)\textbf{I}_n$, so that $P_{\pi}\left(\theta \le 0 \mid \sigma^2 \right)=\frac{1}{2}$. Based on Lemma \ref{inconBF2}, the Bayes factor is then given by
\begin{eqnarray}
\nonumber B_{10} &=& T_{n}\left(-\sqrt{\frac{n^2}{1+(n-1)\rho+t^{-2}(n-1)(1+n+(n-1)\rho)}}\right)^{-1}-1\\
\label{onesidedDim1}& \rightarrow & T_{n}\left(-n(1+(n-1)\rho)^{-\frac{1}{2}}\right)^{-1}-1,
\end{eqnarray}
as $t\rightarrow\infty$, where $T_{\nu}(\cdot)$ denotes the cdf of a univariate Student t distribution with $\nu$ degrees of freedom. Note that as $t\rightarrow-\infty$, $B_{10}$ converges to the reciprocal of \eqref{onesidedDim1}.

\begin{figure}[t]
\centering
\makebox{\includegraphics[width=5.5cm,height=5.1cm]{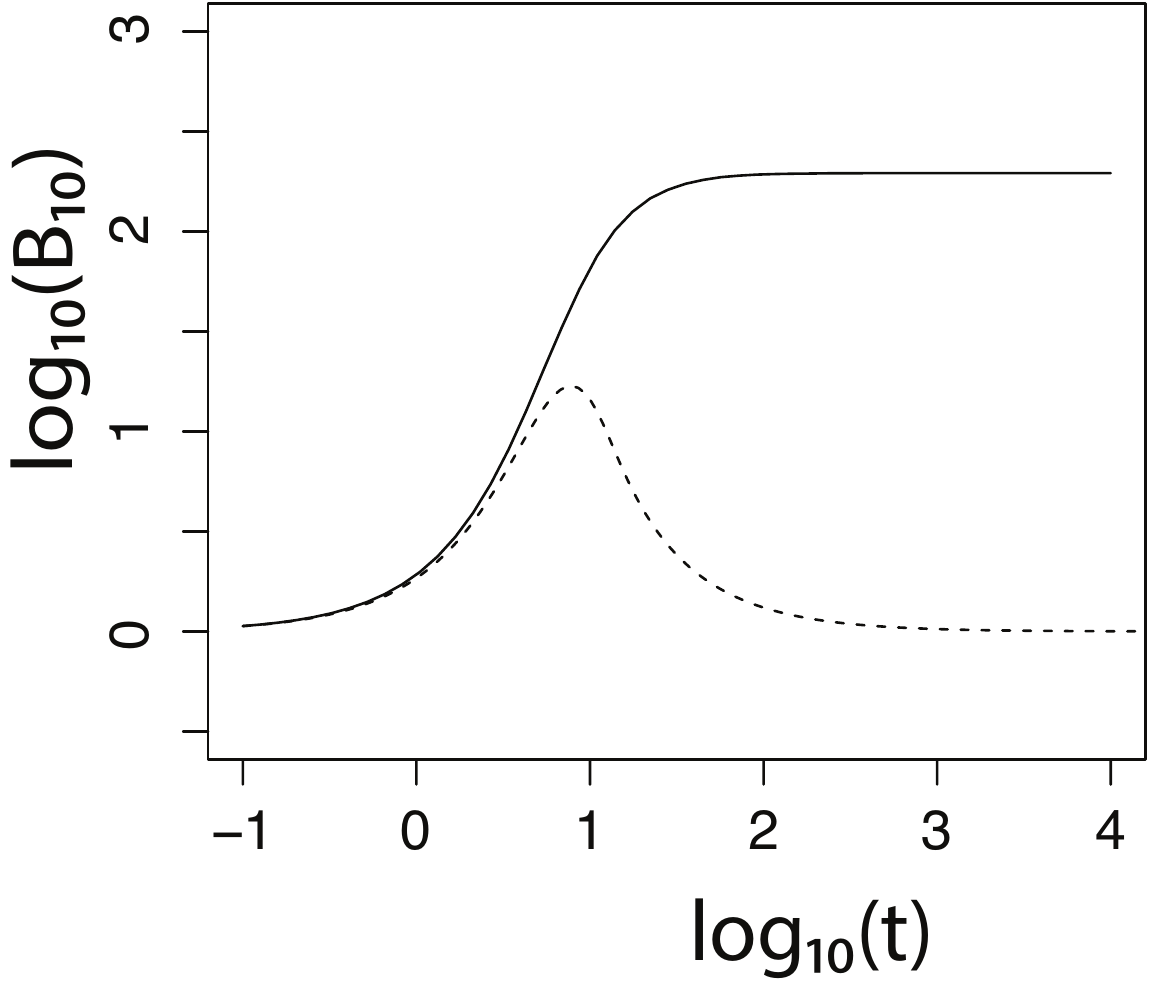}}
\caption{The Bayes factor $B_{10}$ for the one-sided hypothesis test based on the conjugate prior (solid line) and independence prior (dashed line) as a function of $t$-values when $n=7$, $\rho=.5$, $s_{\textbf{y}}^2=n-1=6$, and setting the objective prior to be improper via $\nu=0$.} \label{fig2}
\end{figure}

Table \ref{table2} provides the limiting values of the Bayes factors and Bayes factors in the case of a relatively large t-value of $4$ for different sample sizes and correlations. When comparing Table \ref{table2} with Table \ref{table1}, we can conclude that the practical importance of information inconsistency for one-sided hypothesis testing is considerably less problematic in comparison to the null hypothesis test. Finally Figure \ref{fig2} (solid line) displays the Bayes factor for the one-sided hypothesis test as a function of the t-value based on $n=7$, $\rho=.5$, $s_{\textbf{y}}^2=n-1=6$, and setting the objective improper based on $\nu=0$.

\begin{table}[t]
\caption{Limiting values of the Bayes factor for a one-sided univariate t test as $t\rightarrow\infty$ for different choices of the sample size $n$ and the correlation $\rho$. Additionally Bayes factors and one-sided p-values are given when $t=4$.}
{\small \begin{tabular}{llccccc}
  \hline
& $n$ & 2 & 5 &7 & 10 &20 \\
\hline
\hline
$\rho=0$ & \mbox{limit}  & 9.90 & $486$ & $9.45\times 10^3$ & $1.26\times 10^6$ & $1.85\times 10^{14}$ \\
& $B_{10}$ for $t=4$ & 8.62 & $78.9$ & $199$ & $510$ & $2.40\times 10^3$ \\
  \hline
$\rho=0.5$ & \mbox{limit}  & $7.19$ & $57.2$ & $199$ & $1.21\times10^3$ & $4.02\times 10^5$ \\
& $B_{10}$ for $t=4$ & $6.50$ &  $25.5$ &  $44.7$ & $81.5$ & $238$ \\
  \hline
$\rho \approx 1$ & \mbox{limit}  & $5.83$ & $25.5$ & $59.3$ & $197$ & $8.57\times10^4$ \\
& $B_{10}$ for $t=4$ & 5.37 & 14.7 & 22.4 & 35.2 & 80.9 \\
  \hline
\hline
& one-sided\\
& $p$-value for $t=4$ & 0.078 & 0.008 & 0.0036 & 0.0016 & 0.0038 \\
  \hline
\end{tabular}}
\label{table2}
\end{table}

\subsection{Mixtures of conjugate priors}
We provide the following necessary and sufficient condition for information consistency for a scale mixture of conjugate normal priors in a one-sided hypothesis test.

\begin{lemma} \label{onesided} Let $\bm{\theta} \mid g, \sigma^2 \sim N_{r_1}(\bm{0}_{r_1}, g \sigma^2 \Om )$,  where $\sigma^2$ has the prior specified in (6) and $g$ has a prior with density $\pi(g)$, and let $\bm{w} = E(\bm{\theta} \mid g, \bm{y})$. Assume that if there exists $i$ such that $\widehat{\theta}_i \rightarrow +\infty$, there exists $j$ such that $w_j > 0$. Alternatively, assume that if $\widehat{\theta}_i \rightarrow - \infty$ for all $i$, then $w_i < 0$ for all $i$. [For instance, this condition is satisfied if $\bm{\theta}$ is univariate or  $\bm{\Omega} \propto (\boldsymbol{X}_{\bm\theta}'\boldsymbol{\Sigma}^{-1}\boldsymbol{X}_{\bm\theta})^{-1}$]. Then, the condition
$$
\int_0^\infty (g+1)^{(n-r_1-r_2+\nu)/2} \pi(g) \mathrm{d} g= \infty
$$
is necessary and sufficient for information consistency.
\end{lemma}
\textit{Proof:} See Appendix~\ref{appF}.

\subsection{Independence prior}
The independence semi-conjugate encompassing prior is given by
\begin{eqnarray}
\pi(\bm{\theta},\bm{\gamma},\sigma^2)\propto N_{\boldsymbol{\theta}}(\textbf{0},\bm\Omega)\times\mbox{inv-}\chi^2_{\sigma^2}(s^2,\nu).
\label{indepprior2}
\end{eqnarray}
The truncated priors of $\bm\theta$ under the nonnested hypotheses are as in \eqref{priornonnested}, except that the normalizing constant $P_{\pi}(\bm{\theta}\in\bm{\Theta}_t)$ is the marginal prior probability of $\bm\Theta_t$.

The Bayes factor for the one-sided hypothesis test based on the independence prior can again be expressed as
\begin{equation}
B_{10} = \left(P_{\pi}(\bm{\theta}\le\textbf{0})^{-1}-1\right)^{-1} \left(P_{\pi}(\bm{\theta}\le\textbf{0} \mid \textbf{y})^{-1}-1\right),
\label{B10onesided2}
\end{equation}
but note that the posterior probability is no longer available in closed form.

\begin{lemma}
\label{inconBF3}
As $|\hat{\bm\theta}| \rightarrow \infty$ and at least one coordinate of $\hat{\bm\theta}$ goes to $\infty$, the Bayes factor of $H_1:\bm\theta\not\le\textbf{0}$ versus $H_0:\bm\theta\le\textbf{0}$ based on the independence encompassing prior in \eqref{indepprior2} satisfies
\begin{eqnarray*}
B_{10}&\rightarrow& \left(P_{\pi}(\bm{\theta}\le\textbf{0})^{-1}-1\right)^{-1}.
\end{eqnarray*}
\end{lemma}
\textit{Proof: See Appendix \ref{proof.inconBF3}.}

Thus, as in null hypothesis testing, the independence prior results in a serious violation of information consistency because the evidence in the data of $H_1$ relative to $H_0$ goes to 1 when the evidence against $H_0$ appears to be overwhelming. For completeness, the Bayes factor for the one-sided hypothesis test is also displayed in Figure \ref{fig2} (dashed line), illustrating the extreme form of information inconsistency.

\subsection{Adaptive priors}

An adaptive prior can be specified where the prior covariance matrix of $\bm\theta$ is adapted to the likelihood such that the Bayes factor is maximized for the hypothesis that is supported by the data (i.e., maximize$B_{01}$ if $\hat{\bm\theta}\le\textbf{0}$, and maximize $B_{10}$ elsewhere). Here we show that an adaptive $g$ prior results in an information consistent Bayes factor.


\begin{lemma}\label{lemmaNonnestedAdaptive_g}
The Bayes factor based on the $g$-prior, with $g_{\max}=\arg\max_g \{B_{01}\}$ if $\hat{\bm\theta}\le\textbf{0}$ and $g_{\max}=\arg\max_g \{B_{10}\}$ if $\hat{\bm\theta}\not\le\textbf{0}$, is information consistent for one-sided hypothesis testing.
\end{lemma}
\textit{Proof:} A proof is given in Appendix \ref{proofNonnestedAdaptive_g}.\\

As shown in the proof, the choice for $g$ that maximizes the Bayes factor is obtained by letting $g$ go to $\infty$. This was also proposed by \cite{Mulder:2014a} as a nonadaptive solution. As a result of letting the prior variances go to infinity, the posterior is not shrunk towards the prior mean, which is sufficient to establish information consistency. Therefore the methods of \cite{Mulder:2014b} and \cite{Gu:2014} are also information consistent. A potential issue of letting $g$ go to infinity is that the marginal likelihoods under $H_0$ and $H_1$ go to 0 in the limit. However because the Bayes factor in \eqref{indepprior2} converges in the limit, with posterior probabilities that are computed using flat priors and prior probabilities that are based on the prior covariance structure, the outcome seems a reasonable default quantification of the relative evidence for a one-sided test.

\section{Multiple hypothesis testing}\label{multipleHypo}
We consider the following definition of the (in)formation consistency in a multiple testing problem.

\begin{definition}
\label{incon2}
A Bayes factor is {\em information consistent}, for a multiple hypothesis test of $H_0:\vet{\theta}=\textbf{0}$ versus $H_1:\vet{\theta}\in\vet{\Theta}_1=\{\vet{\theta} \mid \vet{\theta}\le\textbf{0}\mbox{ and } \vet{\theta}\not = \textbf{0}\}$ versus $H_2:\vet{\theta}\in\vet{\Theta}_2=\{\vet{\theta} \mid \vet{\theta}\not \le \textbf{0} \}$, if $B_{20},B_{21}\rightarrow\infty$ as $|\hat{\bm\theta}| \rightarrow \infty$ with at least one coordinate of $\hat{\bm\theta}$ going to $\infty$, and $B_{10},B_{12}\rightarrow \infty$, as all coordinates of $\hat{\bm\theta}$ go to $-\infty$. If this does not hold, the Bayes factor is called information inconsistent.
\end{definition}

As the conjugate and independent semi-conjugate priors resulted in information inconsistent Bayes factors for the one-sided hypothesis test, this automatically implies that these priors result in information inconsistency for the multiple hypothesis test. A specific case that is interesting to mention is when setting $\nu_0>\nu$ when using conjugate priors. This setting of the prior degrees of freedom results in information consistency for the precise hypothesis test. As $|\hat{\bm\theta}|\rightarrow\infty$, and at least one coordinate of $\hat{\bm\theta}$ goes to $\infty$, the Bayes factor $B_{20}$ goes to $\infty$ (a consequence of Lemma \ref{lemma1}), and the Bayes factor $B_{21}$ goes to $B^*_{21}<\infty$ (a consequence of Lemma \ref{inconBF2}). In the case of the simple univariate t-test, this implies that, as $t\rightarrow\infty$, the support for a negative effect, $H_1:\theta<0$, relative to no effect, $H_0:\theta=0$, will go to $\infty$. Thus, even though setting more prior degrees of freedom under the precise hypothesis than under the alternative results in information consistent behavior for the precise test (Lemma \ref{lemma1}), generalizing this prior specification to the multiple hypothesis test results in a dramatic form of information inconsistency.

Finally note that Lemma \ref{genmixt} and \ref{onesided} give the necessary and sufficient conditions for the mixing distribution of the scale mixture of conjugate priors to be information consistent in the multiple testing problem.

\section{Conclusions}\label{concl}
The first major conclusion is that information inconsistency is ubiquitous in hypothesis testing and model selection, when conjugate priors are used. It happens in standard null hypothesis testing and one-sided testing; it happens with proper and improper conjugate priors; and it happens with independence conjugate priors almost always. The practical importance of the problem varies over different situations; it will primarily be a practical problem when the sample is small relative to the number of free parameters and there is high correlation between the observations. But, even in other cases, we consider information inconsistency to be highlighting a logical flaw that might have other serious consequences and is, hence, something to be avoided.

The second major conclusion is that use of either fatter-tailed priors (including appropriate mixtures of g-priors) or adaptive priors typically result in information consistency. This is not as surprising as the almost complete lack of information consistency for conjugate priors, in that previous particular fatter-tailed priors (such as the Zellner-Siow prior) had been shown to be information consistent. Still, the generality in which such priors can be shown to be information consistent is highly comforting.

It should be noted that, when proper priors yield information inconsistency, a logical flaw in Bayesian analysis is not being discovered; if one truly believed the priors were correct, then one should behave in an information inconsistent manner. But one rarely accurately knows features of the priors -- such as their tail behaviors -- that determine information inconsistency. Thus the intuitive appeal of information consistency can be used as a significant aid to selection
of such prior features.

\appendix



\section{Proof of Lemma \ref{genmixt}} \label{appA}

Denote:
\begin{align*}
\widehat{\boldsymbol{\theta}} &= \left(\textbf{X}_{\bm\theta}'\boldsymbol{\Sigma}^{-1}\textbf{X}_{\bm\theta}\right)^{-1}\textbf{X}_{\bm\theta}'\boldsymbol{\Sigma}^{-1}\textbf{y} \\
s^2_{\textbf{y}} &=(\textbf{y}-\textbf{X}_{\bm\theta}\widehat{\boldsymbol{\theta}}-\textbf{X}_{\bm\gamma}\widehat{\boldsymbol{\gamma}})'\bm\Sigma^{-1}(\textbf{y}-\textbf{X}_{\bm\theta}\widehat{\boldsymbol{\theta}}-\textbf{X}_{\bm\gamma}\widehat{\boldsymbol{\gamma}}) \\
\mathsf{SSE}_0 &= s_0^2\nu_0+s^2_{\textbf{y}} \\
\mathsf{SSE}_1 &= s_1^2\nu_1+s^2_{\textbf{y}} \\
\mathsf{SSR} &= \widehat{\boldsymbol{\theta}}'
\textbf{X}_{\bm\theta}'\boldsymbol{\Sigma}^{-1}\textbf{X}_{\bm\theta}\widehat{\boldsymbol{\theta}} \\
\mathcal{I}_\theta &= \mathbf{X}_\theta' \bm{\Sigma}^{-1} \mathbf{X}_\theta  \\
p_0 &= r_2 - \nu_0 \\
p_1 &= r_2 - \nu_1
\end{align*}
Throughout, we use the following notation for functions $a,b$:
\begin{itemize}
    \item $a(g, \hatth) \lesssim b(g, \hatth)$ if and only if there exists $0 < M < \infty$ which doesn't depend on $g$ or $\hatth$ such that $a(g, \hatth) \le M b(g,\hatth).$
    \item $a(g, \hatth) \gtrsim b(g, \hatth)$ if and only if there exists $0 < M < \infty$ which doesn't depend on $g$ or $\hatth$ such that $a(g, \hatth) \ge M b(g,\hatth).$
    \item $a(g, \hatth) \asymp b(g, \hatth)$ if and only if $a(g, \hatth) \lesssim b(g, \hatth)$ and $a(g, \hatth) \gtrsim b(g, \hatth)$.
\end{itemize}
Before we prove Lemma~\ref{genmixt}, we prove an auxiliary result
\begin{lemma}  \label{boundh}
Let
$$
h(g) = |g \Om + \J|^{-1/2} [\mathsf{SSE}_1 + \hatth'(g\Om + \J)^{-1} \hatth ]^{-(n-p_1)/2},
$$
then, there exist $0 < d_l < d_u < \infty$ such that
$$
 \frac{(g+d_l)^{(n-p_1-r_1)/2}}{[(g+d_l) \mathsf{SSE}_1+ \hatth' \Om^{-1} \hatth ]^{(n-p_1)/2}} \lesssim h(g) \lesssim \frac{(g+d_u)^{(n-p_1-r_1)/2}}{[(g+d_u)\mathsf{SSE}_1+\hatth' \Om^{-1} \hatth ]^{(n-p_1)/2}}
$$
\end{lemma}
\textit{Proof:} Consider the matrix factorization
$$
\mathcal{I}_\theta^{-1} + g \bm{\Omega} = \bm{\Omega}^{1/2} [ \bm{\Omega}^{-1/2} \mathcal{I}^{-1}_\theta \bm{\Omega}^{-1/2} + g I_{r_1}] \bm{\Omega}^{1/2},
$$
and take the eigendecomposition $\bm{\Omega}^{-1/2} \mathcal{I}^{-1}_\theta \bm{\Omega}^{-1/2} = \bm{O} \bm{D} \bm{O}'$, where $\bm{O}$ is orthogonal and $\bm{D}$ diagonal with elements $0 < d_l < d_i < d_u < \infty$. Then, we can rewrite
$$
\mathcal{I}_\theta^{-1} + g \bm{\Omega} = \Om^{1/2} \bm{O} [ \bm{D} + g I_{r_1}] \bm{O}' \bm{\Omega}^{1/2}.
$$
We can bound
$$
\hatth' \Om^{-1} \hatth/(d_u+g)  \le \hatth'(g\Om + \J)^{-1} \hatth \le \hatth' \Om^{-1} \hatth/(d_l+g)
$$
and
$$
| g \Om + \J|^{-1/2} \propto | \bm{D} + g I_{r_1}|^{-1/2} \in [(g+d_u)^{-r_1/2}, (g+d_l)^{-r_1/2}],
$$
so
$$
h(g) \lesssim \frac{(d_u+g)^{(n-p_1)/2} (d_l+g)^{-r_1/2}}{[(d_u+g)\mathsf{SSE}_1+\hatth' \Om^{-1} \hatth ]^{(n-p_1)/2}} \lesssim \frac{(d_u+g)^{(n-p_1-r_1)/2}}{[(d_u+g)\mathsf{SSE}_1+\hatth' \Om^{-1} \hatth ]^{(n-p_1)/2}}.
$$
Similarly, we can find the lower bound
$$
h(g) \gtrsim \frac{(g+d_l)^{(n-p_1-r_1)/2}}{[(g+d_l) \mathsf{SSE}_1+ \hatth' \Om^{-1} \hatth ]^{(n-p_1)/2}}.
$$
Now, we prove Lemma~\ref{genmixt} arguing by cases.  \\
 \underline{Case $\nu_0 > \nu_1$}
 Applying the lower bound in Lemma~\ref{boundh},
 \begin{align*}
 B_{10} \gtrsim   \frac{\left[\mathsf{SSE}_0+ \mathsf{SSR} \right]^{(n-p_0)/2}}{(\hatth' \Om^{-1} \hatth)^{(n-p_1)/2} } \int_0^\infty \frac{(g+d_l)^{(n-p_1-r_1)/2}}{[(g+d_l) \frac{\mathsf{SSE}_1}{\hatth' \Om^{-1} \hatth}+ 1 ]^{(n-p_1)/2}} \pi(\mathrm{d} g).
\end{align*}
Since $p_0 < p_1,$ the term outside the integral goes to infinity as $\normth^2 \rightarrow \infty$, and
 by Fatou's lemma,
 $$
 \lim \inf_{\normth^2 \rightarrow \infty} \int_0^\infty \frac{(g+d_l)^{(n-p_1-r_1)/2}}{[(g+d_l) \frac{\mathsf{SSE}_1}{\hatth' \Om^{-1} \hatth}+ 1 ]^{(n-p_1)/2}} \pi(\mathrm{d} g) \ge \int_0^\infty (g+d_l)^{(n-p_1-r_1)/2} \pi( \mathrm{d} g),
 $$
 which is clearly bounded away from 0 for any prior on $g$ with positive support, so any such prior yields an information-consistent $B_{10}$ whenever $\nu_0 > \nu_1.$ \\

 \noindent \underline{Case $\nu_0 = \nu_1$} Applying the lower bound in Lemma~\ref{boundh} and Fatou's lemma as we did for the case $\nu_0 > \nu_1$:
 $$
 \lim_{\normth^2 \rightarrow \infty} B_{10} \gtrsim  \int_0^\infty (g+d_l)^{(n-p_1-r_1)/2} \pi( \mathrm{d} g) \lim_{\normth^2 \rightarrow \infty}  \frac{\left[\mathsf{SSE}_0+ \mathsf{SSR} \right]^{(n-p_0)/2}}{[\hatth' \Om^{-1} \hatth]^{(n-p_1)/2} } .
 $$
The limit is $O(1)$, so a sufficient condition for information consistency is
$$
\int_0^\infty (g+d_l)^{(n-p_1-r_1)/2} \pi( \mathrm{d} g) \asymp \int_0^\infty (g+1)^{(n-p_1-r_1)/2} \pi( \mathrm{d} g)  = \infty,
$$
as required. \\

\noindent \underline{Case $\nu_0 < \nu_1$} In this case, we apply the upper bound in Lemma~\ref{boundh}:
 $$
 B_{10} \lesssim \frac{\left[\mathsf{SSE}_0+ \mathsf{SSR} \right]^{(n-p_0)/2}}{(\hatth' \Om^{-1} \hatth)^{(n-p_1)/2} } \int_0^\infty \frac{(g+d_u)^{(n-p_1-r_1)/2}}{[(g+d_u) \frac{\mathsf{SSE}_1}{\hatth' \Om^{-1} \hatth}+ 1 ]^{(n-p_1)/2}} \pi(\mathrm{d} g).
 $$
 The term outside the integral goes to 0, so a necessary condition for information consistency is that the integral be infinite. We can bound the integral:
 $$
 \int_0^\infty \frac{(g+d_u)^{(n-p_1-r_1)/2}}{[(g+d_u) \frac{\mathsf{SSE}_1}{\hatth' \Om^{-1} \hatth}+ 1 ]^{(n-p_1)/2}} \pi(\mathrm{d} g) \le  \int_0^\infty (g+d_u)^{(n-p_1-r_1)/2} \, \pi(\mathrm{d} g),
 $$
 so a necessary condition for information consistency is
 $$
 \int_0^\infty (g+d_u)^{(n-p_1-r_1)/2} \, \pi(\mathrm{d} g) \asymp \int_0^\infty (g+1)^{(n-p_1-r_1)/2} \pi( \mathrm{d} g)  = \infty,
$$
as required.

\section{Proof of Lemma \ref{tail}} \label{appB}

Throughout, we use the notation in Appendix~\ref{appA}.

\noindent \underline{Case 1.}  Suppose there exists $M < \infty$ such that for all $g \ge M$, $\pi(g) \gtrsim g^{-\alpha}$ for $\alpha >1$ and $p_0 > p_1$. Then, we apply the lower bound in Lemma~\ref{boundh}:
$$
B_{10} \gtrsim   \left[\mathsf{SSE}_0+ \mathsf{SSR} \right]^{(n-p_0)/2} \int_M^\infty \frac{(g+d_l)^{(n-p_1-r_1)/2-\alpha}}{[(g+d_l) \mathsf{SSE}_1 + \hatth' \Om^{-1} \hatth]^{(n-p_1)/2}} \mathrm{d} g.
$$
Now, note that for any $K,d > 0$ with $1-d < K$,
 $$
0 \le  \lim_{\normth^2 \rightarrow \infty} \int_{\min(0,1-d)}^K \frac{(g+d)^{(n-p_1-r_1)/2-\alpha}}{[(g+d) \mathsf{SSE}_1 + \hatth' \Om^{-1} \hatth]^{(n-p_1)/2}} \mathrm{d} g \lesssim \lim_{\normth^2 \rightarrow \infty}[\hatth' \Om^{-1} \hatth]^{-(n-p_1)/2} = 0,
 $$
so
$$
\lim_{\normth^2 \rightarrow \infty} \int_M^\infty \frac{(g+d_l)^{(n-p_1-r_1)/2-\alpha}}{[(g+d_l) \mathsf{SSE}_1 + \hatth' \Om^{-1} \hatth]^{(n-p_1)/2}} \mathrm{d} g =  \lim_{\normth^2 \rightarrow \infty} \int_{1-d_l}^\infty \frac{(g+d_l)^{(n-p_1-r_1)/2-\alpha}}{[(g+d_l) \mathsf{SSE}_1 + \hatth' \Om^{-1} \hatth]^{(n-p_1)/2}} \mathrm{d} g.
$$
Plugging in:
 \begin{align*}
 \lim_{\normth^2 \rightarrow \infty}  B_{10} &\gtrsim  \lim_{\normth^2 \rightarrow \infty}   \left[\mathsf{SSE}_0+ \mathsf{SSR} \right]^{(n-p_0)/2} \int_{1-d_l}^\infty \frac{(g+d_l)^{(n-p_1-r_1)/2-\alpha}}{[(g+d_l) \mathsf{SSE}_1 + \hatth' \Om^{-1} \hatth]^{(n-p_1)/2}} \mathrm{d} g \\
 &\propto  \lim_{\normth^2 \rightarrow \infty} \frac{(\mathsf{SSE}_0 + \mathsf{SSR})^{(n-p_0)/2}}{\mathsf{SSE}_1^{(n-p_1)/2}} {}_{2}F_1\left(\frac{n-p_1}{2},\frac{r_1}{2}+\alpha-1; \frac{r_1}{2}+\alpha; \frac{{- \hatth' \Om^{-1} \hatth}}{\mathsf{SSE}_1} \right).
\end{align*}
Using the identity
$$
{}_2F_1 (a,b;c;z) = (1-z)^{-b} {}_2F_1 \left (b,c-a;c;\tfrac{z}{z-1} \right ),
$$
we have
\begin{align*}
\lim_{\normth^2 \rightarrow \infty} B_{10} &\gtrsim \lim_{\normth^2 \rightarrow \infty} \frac{(\mathsf{SSE}_0 + \mathsf{SSR})^{(n-p_0)/2} {}_{2}F_1\left(\frac{r_1}{2}+\alpha-1,\frac{r_1-(n-p_1)}{2}+\alpha;\frac{r_1}{2}+\alpha ;  R^2 \right)}{\mathsf{SSE}_1^{(n-p_1)/2} \left[ 1 + \frac{ \hatth' \Om^{-1} \hatth}{\mathsf{SSE}_1} \right]^{(r_1/2)+\alpha-1}} ,
\end{align*}
where $R^2 =  \hatth' \Om^{-1} \hatth/( \hatth' \Om^{-1} \hatth+\mathsf{SSE}_1) \rightarrow 1$ as $\normth^2 \rightarrow \infty$. If $\alpha < (n-p_1-r_1)/2+1$ (which is satisfied because $\alpha < (n-p_0-r_1)/2$ and $p_0 > p_1$ by assumption), the limit of the hypergeometric function as $R^2 \rightarrow 1$ is a constant (by Gauss' theorem). From here, it is immediate to conclude that $B_{10}$ is information consistent whenever the lower bound is infinite, which occurs for $\alpha < (n-p_0-r_1)/2+1,$ as required. \\

\noindent \underline{Case 2.} Suppose there exists $M' < \infty$ such that for all $g \ge M'$, $\pi(g) \lesssim g^{-\alpha}$ for $\alpha >1$ and $p_0 > p_1$. Then, by Lemma~\ref{boundh}:
$$
B_{10} \lesssim \left[\mathsf{SSE}_0+ \mathsf{SSR} \right]^{(n-p_0)/2} \int_M^\infty \frac{(g+d_u)^{(n-p_1-r_1)/2-\alpha}}{[(g+d_u) \mathsf{SSE}_1 + \hatth' \Om^{-1} \hatth]^{(n-p_1)/2}} \mathrm{d} g.
$$
As argued in Case 1,
\begin{align*}
\lim_{\normth^2 \rightarrow \infty} \int_M^\infty \frac{(g+d_u)^{(n-p_1-r_1)/2-\alpha}}{[(g+d_u) \mathsf{SSE}_1 + \hatth' \Om^{-1} \hatth]^{(n-p_1)/2}} \mathrm{d} g = \lim_{\normth^2 \rightarrow \infty} \int_{1-d_u}^\infty \frac{(g+d_u)^{(n-p_1-r_1)/2-\alpha}}{[(g+d_u) \mathsf{SSE}_1 + \hatth' \Om^{-1} \hatth]^{(n-p_1)/2}} \mathrm{d} g,
\end{align*}
and carrying out the same computations as in Case 1:
\begin{align*}
\lim_{\normth^2 \rightarrow \infty} B_{10} &\lesssim \lim_{\normth^2 \rightarrow \infty} \frac{(\mathsf{SSE}_0 + \mathsf{SSR})^{(n-p_0)/2} {}_{2}F_1\left(\frac{r_1}{2}+\alpha-1,\frac{r_1-(n-p_1)}{2}+\alpha;\frac{r_1}{2}+\alpha ;  R^2 \right)}{\mathsf{SSE}_1^{(n-p_1)/2} \left[ 1 + \frac{ \hatth' \Om^{-1} \hatth}{\mathsf{SSE}_1} \right]^{(r_1/2)+\alpha-1}}.
\end{align*}
If $(n-p_0-r_1)/2 +1 \le \alpha < (n-p_1-r_1)/2+1$, the limit of the hypergeometric function is $O(1)$ and $B_{10}$ is information inconsistent. If $\alpha \ge (n-p_1-r_1)/2+1$, the necessary condition of Lemma~\ref{genmixt} implies that $B_{10}$ is information inconsistent. Therefore, $B_{10}$ is information inconsistent whenever $\alpha \ge (n-p_0-r_1)/2+1,$ as required.

\section{Proof of Lemma \ref{lemma.indconj}}\label{proof.indconj}
Break the integral in $B_{10}$ into the two regions $R_1=\{\boldsymbol{\theta}: \, |\boldsymbol{\theta}|^2 \leq |\hat{\boldsymbol{\theta}}|\}$
and $R_2 = \{\boldsymbol{\theta}: \, |\boldsymbol{\theta}|^2 > |\hat{\boldsymbol{\theta}}|\}$. It is easy to see that, for any fixed $\epsilon >0$, there is a $K_{\epsilon}$ such that, for $|\hat{\boldsymbol{\theta}}| > K_{\epsilon}$ and $\boldsymbol{\theta} \in R_1$,
\begin{eqnarray*}
(1-\epsilon)\left(\hat{\boldsymbol{\theta}}'\textbf{X}'_{\bm\theta}\boldsymbol{\Sigma}^{-1}
\textbf{X}_{\bm\theta}\hat{\boldsymbol{\theta}} \right)^{-\frac{n-r_2+\nu_1}{2}}
&<& \left(\nu_1 s_1^2+s^2_{\textbf{y}}+(\boldsymbol{\theta}-\hat{\boldsymbol{\theta}})'\textbf{X}'_{\bm\theta}\boldsymbol{\Sigma}^{-1}
\textbf{X}_{\bm\theta}(\boldsymbol{\theta}-\hat{\boldsymbol{\theta}}) \right)^{-\frac{n-r_2+\nu_1}{2}} \\
&<& (1+\epsilon)\left(\hat{\boldsymbol{\theta}}'\textbf{X}'_{\bm\theta}\boldsymbol{\Sigma}^{-1}
\textbf{X}_{\bm\theta}\hat{\boldsymbol{\theta}} \right)^{-\frac{n-r_2+\nu_1}{2}} \,.
\end{eqnarray*}
Thus, letting $P(R_1)$ denote the probability of $R_1$ under the $N_{\boldsymbol{\theta}}(\textbf{0},\bm\Omega)$ density, it follows that,
for $|\hat{\boldsymbol{\theta}}| > K_{\epsilon}$,
\begin{eqnarray*}
&&(1-\epsilon)\left(\hat{\boldsymbol{\theta}}'\textbf{X}'_{\bm\theta}\boldsymbol{\Sigma}^{-1}P(
\textbf{X}_{\bm\theta}\hat{\boldsymbol{\theta}} \right)^{-\frac{n-r_2+\nu_1}{2}} P(R_1) \\
&<& \int_{R_1} \left(\nu_1 s_1^2+s^2_{\textbf{y}}+(\boldsymbol{\theta}-\hat{\boldsymbol{\theta}})'\textbf{X}'_{\bm\theta}\boldsymbol{\Sigma}^{-1}
\textbf{X}_{\bm\theta}(\boldsymbol{\theta}-\hat{\boldsymbol{\theta}}) \right)^{-\frac{n-r_2+\nu_1}{2}} N_{\boldsymbol{\theta}}(\textbf{0},\bm\Omega) d\boldsymbol{\theta}\\
&<& (1+\epsilon)\left(\hat{\boldsymbol{\theta}}'\textbf{X}'_{\bm\theta}\boldsymbol{\Sigma}^{-1}
\textbf{X}_{\bm\theta}\hat{\boldsymbol{\theta}} \right)^{-\frac{n-r_2+\nu_1}{2}} P(R_1)\,.
\end{eqnarray*}
As $|\hat{\boldsymbol{\theta}}| \rightarrow \infty$, the integral over $R_2$ is clearly going to zero exponentially fast, while
$P(R_1) \rightarrow 1$. Since $\epsilon$ can be chosen arbitrarily small, it follows that, as $|\hat{\boldsymbol{\theta}}| \rightarrow \infty$,
$$ \frac {\int \left(\nu_1 s_1^2+s^2_{\textbf{y}}+(\boldsymbol{\theta}-\hat{\boldsymbol{\theta}})'\textbf{X}'_{\bm\theta}\boldsymbol{\Sigma}^{-1}
\textbf{X}_{\bm\theta}(\boldsymbol{\theta}-\hat{\boldsymbol{\theta}}) \right)^{-\frac{n-r_2+\nu_1}{2}} N_{\boldsymbol{\theta}}(\textbf{0},\bm\Omega) d\boldsymbol{\theta}}
{\left(\hat{\boldsymbol{\theta}}'\textbf{X}'_{\bm\theta}\boldsymbol{\Sigma}^{-1}
\textbf{X}_{\bm\theta}\hat{\boldsymbol{\theta}} \right)^{-\frac{n-r_2+\nu_1}{2}} }
\rightarrow 1 \,.
$$
Thus, as $|\hat{\boldsymbol{\theta}}| \rightarrow \infty$,
$$ B_{10} \rightarrow \lim_{|\hat{\boldsymbol{\theta}}| \rightarrow \infty} \frac{C_2 \left(\hat{\boldsymbol{\theta}}'\textbf{X}'_{\bm\theta}\boldsymbol{\Sigma}^{-1}
\textbf{X}_{\bm\theta}\hat{\boldsymbol{\theta}} \right)^{-\frac{n-r_2+\nu_1}{2}}} {{\left(\nu_0 s_0^2+s^2_{\textbf{y}}+\hat{\boldsymbol{\theta}}'\textbf{X}'_{\bm\theta}\boldsymbol{\Sigma}^{-1}\textbf{X}_{\bm\theta}\hat{\boldsymbol{\theta}} \right)^{-\frac{n-r_2+\nu_0}{2}}}} \,,
$$
from which the results stated in the lemma follow directly.

\section{Proof of Lemma~\ref{adaptive}} \label{appD}

\noindent Using the notation in Appendix~\ref{appA} and applying Lemma~\ref{boundh}:
$$
B_{10} \gtrsim  \frac{ (\mathsf{SSE}_0 + \mathsf{SSR})^{(n-p_0)/2}}{[\hatth' \Om^{-1} \hatth]^{n-p_1/2} }  \frac{(g+d_l)^{(n-p_1-r_1)/2}}{ [(g+d_l) \mathsf{SSE}_1/\hatth' \Om^{-1} \hatth +1 ]^{(n-p_1)/2}}
$$
For $g >0,$ the right-hand side is maximized at $\widehat{g} = \max(0,(n-p_1-r_1) \hatth' \Om^{-1} \hatth/(r_1 \mathsf{SSE})-d_l)$. Then,
\begin{align*}
\lim_{\normth^2 \rightarrow \infty} \max_{g \ge 0} B_{10} &\gtrsim  \frac{ (\mathsf{SSE}_0 + \mathsf{SSR})^{(n-p_0)/2}}{[\hatth' \Om^{-1} \hatth]^{(n-p_1)/2} }  \frac{(\widehat{g}+d_l)^{(n-p_1-r_1)/2}}{ [(\widehat{g}+d_l) \mathsf{SSE}_1/\hatth' \Om^{-1} \hatth +1 ]^{(n-p_1)/2}} \\
&\propto \lim_{\normth^2 \rightarrow \infty} \frac{ (\mathsf{SSE}_0 + \mathsf{SSR})^{(n-p_0)/2}}{[\hatth' \Om^{-1} \hatth]^{r_1/2} \mathsf{SSE}^{(n-p_1-r_1)/2} } \\
&= \infty,
\end{align*}
so the adaptive prior is information consistent.

\section{Proof of Lemma \ref{inconBF2}}\label{proof.inconBF2}
The marginal posterior of $\bm\theta$ in the joint space has a multivariate Student t distribution with mean $(\textbf{X}_{\bm\theta}'\bm{\Sigma}^{-1}\textbf{X}_{\bm\theta}+
\bm\Omega^{-1})^{-1}\textbf{X}_{\bm\theta}'\bm{\Sigma}^{-1}\textbf{X}_{\bm\theta}\hat{\bm{\theta}}$, scale matrix $(n+\nu-r_2)^{-1}(s^2\nu+s_{\textbf{y}}^2+\hat{\bm{\theta}}'((\textbf{X}_{\bm\theta}'\bm{\Sigma}^{-1}\textbf{X}_{\bm\theta})^{-1}+\bm\Omega)^{-1}\hat{\bm{\theta}})
(\textbf{X}_{\bm\theta}'\bm{\Sigma}^{-1}\textbf{X}_{\bm\theta}+\bm\Omega^{-1})^{-1}$, and $n+\nu-r_2$ degrees of freedom. Change variables to
$$\bm{\xi}=(n+\nu-r_2)^{1/2}(s^2\nu+s_{\textbf{y}}^2+\hat{\bm{\theta}}'((\textbf{X}_{\bm\theta}'\bm{\Sigma}^{-1}\textbf{X}_{\bm\theta})^{-1}+
\bm\Omega)^{-1}\hat{\bm{\theta}})^{-1/2} \bm{\theta} \,,$$
which has a multivariate Student t distribution with mean
$$ \bm{\xi}^* =\frac{(\textbf{X}_{\bm\theta}'\bm{\Sigma}^{-1}\textbf{X}_{\bm\theta}+
\bm\Omega^{-1})^{-1}\textbf{X}_{\bm\theta}'\bm{\Sigma}^{-1}\textbf{X}_{\bm\theta}\hat{\bm{\theta}}}
{(n+\nu-r_2)^{-1/2}(s^2\nu+s_{\textbf{y}}^2+\hat{\bm{\theta}}'((\textbf{X}_{\bm\theta}'\bm{\Sigma}^{-1}\textbf{X}_{\bm\theta})^{-1}+
\bm\Omega)^{-1}\hat{\bm{\theta}})^{1/2}} \,,$$
scale matrix $(\textbf{X}_{\bm\theta}'\bm{\Sigma}^{-1}\textbf{X}_{\bm\theta}+\bm\Omega^{-1})^{-1}$, and  $n+\nu-r_2$ degrees of freedom.
Note that
$$P_{\pi}(\bm{\theta}\le\textbf{0} \mid \textbf{y}) = P_{\pi}(\bm{\xi}\le\textbf{0} \mid \textbf{y}) \,.$$
It is easy to see that $\bm{\xi}^*$
lies in a fixed compact set $C$ for any $\hat{\bm{\theta}}$, from which it is immediate that $P_{\pi}(\bm{\xi}\le\textbf{0} \mid \textbf{y})$ is bounded away from 0 and 1.

The second part of the lemma follows immediately from letting $c \rightarrow \infty$ in the expression for $ \bm{\xi}^*$.

\section{Proof of Lemma \ref{onesided}} \label{appF}

\noindent Throughout, we use the notation in Appendix~\ref{appA}.

\noindent \underline{Sufficient condition:} \\
We start with the case where there exists $\widehat{\theta}_i \rightarrow +\infty$; we treat the case where all $\widehat{\theta}_i \rightarrow -\infty$ later. \\

\noindent We can write:
\begin{align*}
\lim_{\lVert \widehat{\bm{\theta}} \rVert^2 \rightarrow \infty} P( \bm{\theta} \le \bm{0} \mid \bm{y} ) &= \lim_{\lVert \widehat{\bm{\theta}} \rVert^2 \rightarrow \infty} \int_0^\infty P( \bm{\theta} \le \bm{0} \mid g, \bm{y} )  p(g \mid \bm{y}) \mathrm{d} g \\
&= \lim_{\lVert \widehat{\bm{\theta}} \rVert^2 \rightarrow \infty} \frac{1}{p(\bm{y})} \int_0^\infty P( \bm{\theta} \le \bm{0} \mid g, \bm{y} )  p( \bm{y} \mid g) \pi(\mathrm{d} g) \\
&\propto\lim_{\lVert \widehat{\bm{\theta}} \rVert^2 \rightarrow \infty}  \frac{1}{p(\bm{y})}  \int_0^\infty  P( \bm{\theta} \le \bm{0} \mid g, \bm{y} )  h(g) \, \pi(\mathrm{d} g),
\end{align*}
with $h$ as defined in Lemma~\ref{boundh} (but noting that, in this case, the notation is $\nu_1 = \nu$). Letting $p = \nu-r_2$ and using the upper bound in Lemma~\ref{boundh},
\begin{align*}
\lim_{\lVert \widehat{\bm{\theta}} \rVert^2 \rightarrow \infty} P( \bm{\theta} \le \bm{0} \mid \bm{y} )& \lesssim \lim_{\lVert \widehat{\bm{\theta}} \rVert^2 \rightarrow \infty} \frac{[\hatth'\Om^{-1} \hatth]^{-(n-p)/2}}{p(\bm{y})} \int_0^\infty \frac{ P( \bm{\theta} \le \bm{0} \mid g, \bm{y} ) \, (g+d_u)^{(n-p-r_1)/2}}{[(g+d_u) \mathsf{SSE}_1/\hatth'\Om^{-1} \hatth+ 1 ]^{(n-p)/2}}  \, \pi(\mathrm{d} g)
\end{align*}
From Lemma~\ref{inconBF2}, we know that
$$
P(\th \le \bm{0} \mid g, \bm{y}) = P(\bm{\xi} \le \bm{0} \mid g, \bm{y}),
$$
where $\bm{\xi}$ has a multivariate Student-$t$ distribution, with location and scale
\begin{align*}
\bm{w} &= (\I + \Om^{-1}/g)^{-1} \I \hatth \\
\bm{m} &= \frac{(n+\nu-r_2)^{1/2} \, \bm{w} }{[\mathsf{SSE}_1 + \hatth' (\J + g \Om)^{-1} \hatth ]^{1/2}} \,  \\
\bm{S} &= (\I + \Om^{-1}/g)^{-1}.
\end{align*}
We factor
\begin{align*}
\bm{S}  = \Om^{1/2} (\Om^{1/2} \I \Om^{1/2} + I_{r_1}/g)^{-1} \Om^{1/2} = \Om^{1/2} \bm{O}' (\bm{D}^{-1} + I_{r_1}/g)^{-1} \bm{O} \Om^{1/2},
\end{align*}
where $\bm{O}$ is orthogonal and $\bm{D}$ is diagonal (with positive entries) as defined in Lemma~\ref{boundh}.
Therefore, for a fixed coordinate $j$,
$$
\bm{S}_{jj} \in \left[ \frac{g}{g/d_{l} + 1} \Om_{jj}, \frac{g}{g/d_{u} + 1} \Om_{jj} \right],
$$
so $0  < \bm{S}_{jj} < \infty$ for $g > 0$. Using the same factorizations, we obtain  $\lVert \bm{w} \rVert^{2} \propto \hatth' \Om \Om \hatth$ for $g > 0$. Plugging this in and factorizing the denominator in $\bm{m}$ in a similar manner, we obtain
\begin{align*}
\bm{m} &= \frac{(n+\nu-r_2)^{1/2} \,  \lVert \bm{w} \rVert }{[\mathsf{SSE}_1 + \hatth' (\J + g\Om)^{-1} \hatth ]^{1/2}} \frac{\bm{w}}{\lVert \bm{w} \rVert } \\
&\propto \frac{ \,    (\hatth' \Om \Om \hatth)^{1/2} }{[\mathsf{SSE}_1 + \hatth' (\J + g\Om)^{-1} \hatth ]^{1/2}} \frac{\bm{w}}{\lVert \bm{w} \rVert }.
\end{align*}
If we choose a coordinate $j$ such that $w_j > 0$ (which exists by assumption), using the lower bound in Lemma~\ref{boundh},
$$
\bm{m}_j \gtrsim \frac{  (g+d_l)^{1/2} (\hatth' \Om \Om \hatth)^{1/2} }{\left[ (g+d_l) \mathsf{SSE}_1 + \hatth \Om^{-1} \hatth  \right]^{1/2}} \gtrsim \frac{(g+d_l)^{1/2}}{\left[(g+d_l) \mathsf{SSE}_1/\hatth'\Om^{-1} \hatth + 1 \right]^{1/2}}
$$
Now, 
\begin{align*}
\lim_{\lVert \widehat{\bm{\theta}} \rVert^2 \rightarrow \infty} P( \bm{\theta} \le \bm{0} \mid \bm{y} )& \lesssim \lim_{\lVert \widehat{\bm{\theta}} \rVert^2 \rightarrow \infty} \frac{[\hatth'\Om^{-1} \hatth]^{-(n-p)/2}}{p(\bm{y})} \int_0^\infty \frac{ P( T_{n-p} \ge \bm{m}_j / \sqrt{\bm{S}_{jj}} ) \, (g+d_u)^{(n-p-r_1)/2}}{[(g+d_u) \mathsf{SSE}_1/\hatth'\Om^{-1} \hatth+ 1 ]^{(n-p)/2}}  \, \pi(\mathrm{d} g)
\end{align*}
where $T_{n-p}$ is a central Student-$t$ with $n-p$ degrees of freedom. Let $\varepsilon > 0,$ then 
$$
\int_0^\varepsilon \frac{ P( T_{n-p} \ge \bm{m}_j / \sqrt{\bm{S}_{jj}} ) \, (g+d_u)^{(n-p-r_1)/2}}{[(g+d_u) \mathsf{SSE}_1/\hatth'\Om^{-1} \hatth+ 1 ]^{(n-p)/2}}  \, \pi(\mathrm{d} g) \le (\varepsilon+d_u)^{(n-p-r_1)/2},
$$
so 
\begin{align*}
\lim_{\lVert \widehat{\bm{\theta}} \rVert^2 \rightarrow \infty} P( \bm{\theta} \le \bm{0} \mid \bm{y} )& \lesssim \lim_{\lVert \widehat{\bm{\theta}} \rVert^2 \rightarrow \infty} \frac{[\hatth'\Om^{-1} \hatth]^{-(n-p)/2}}{p(\bm{y})} \int_\varepsilon^\infty \frac{ P( T_{n-p} \ge \bm{m}_j / \sqrt{\bm{S}_{jj}} ) \, (g+d_u)^{(n-p-r_1)/2}}{[(g+d_u) \mathsf{SSE}_1/\hatth'\Om^{-1} \hatth+ 1 ]^{(n-p)/2}}  \, \pi(\mathrm{d} g).
\end{align*}
Therefore, we can plug in our bounds for $\bm{m}_j$ and $\bm{S}_{jj}$, which are bounded away from 0 whenever $g > 0$. Using the tail bound
$$
P(T_{n-p} \ge x ) \lesssim \frac{1}{x(1 + x^2/\nu)^{(n-p-1)/2}}  \lesssim x^{-(n-p)}
$$
and our previous work, we obtain
\begin{align*}
\lim_{\lVert \widehat{\bm{\theta}} \rVert^2 \rightarrow \infty} P( \bm{\theta} \le \bm{0} \mid \bm{y} ) &\lesssim \lim_{\lVert \widehat{\bm{\theta}} \rVert^2 \rightarrow \infty}  \frac{[\hatth'\Om^{-1} \hatth]^{-(n-p)/2}}{p(\bm{y})} \int_\varepsilon^\infty (g+d_u)^{-r_1/2} \, \pi(\mathrm{d}g) \\ &\propto \lim_{\lVert \widehat{\bm{\theta}} \rVert^2 \rightarrow \infty}  \frac{[\hatth'\Om^{-1} \hatth]^{-(n-p)/2}}{p(\bm{y})}.
\end{align*}
Clearly
$$
\lim_{\lVert \widehat{\bm{\theta}} \rVert^2 \rightarrow \infty}  \frac{[\hatth'\Om^{-1} \hatth]^{-(n-p)/2}}{p(\bm{y})} = 0 \Leftrightarrow \lim_{\lVert \widehat{\bm{\theta}} \rVert^2 \rightarrow \infty} [\hatth'\Om^{-1} \hatth]^{(n-p)/2} \, p(\bm{y}) = \infty
$$
and
\begin{align*}
 \lim_{\lVert \widehat{\bm{\theta}} \rVert^2 \rightarrow \infty} [\hatth'\Om^{-1} \hatth]^{(n-p)/2} \, p(\bm{y}) &=  \lim_{\lVert \widehat{\bm{\theta}} \rVert^2 \rightarrow \infty} [\hatth'\Om^{-1} \hatth]^{(n-p)/2} \, \int_0^\infty p(\bm{y} \mid g) \, \pi(\mathrm{d} g) \\
 &\propto \lim_{\lVert \widehat{\bm{\theta}} \rVert^2 \rightarrow \infty} [\hatth'\Om^{-1} \hatth]^{(n-p)/2} \, \int_0^\infty h(g) \, \pi(\mathrm{d} g) \\
 &\gtrsim \lim_{\lVert \widehat{\bm{\theta}} \rVert^2 \rightarrow \infty} \int_0^{\infty}  \frac{(g+d_l)^{(n-p-r_1)/2}}{[(g+d_l) \mathsf{SSE}_1/\hatth' \Om^{-1} \hatth+ 1]^{(n-p)/2}} \pi(\mathrm{d}g) \\
&\gtrsim  \int_0^{\infty}  \lim \inf_{\lVert \widehat{\bm{\theta}} \rVert^2 \rightarrow \infty} \frac{(g+d_l)^{(n-p-r_1)/2}}{[(g+d_l) \mathsf{SSE}_1/\hatth' \Om^{-1} \hatth+ 1]^{(n-p)/2}} \pi(\mathrm{d}g) \\
&= \int_0^{\infty}  (g+d_l)^{(n-p-r_1)/2} \pi(\mathrm{d}g) \\
&\asymp  \int_0^{\infty}  (g+1)^{(n-p-r_1)/2} \pi(\mathrm{d}g).
\end{align*}
Therefore, if the integral above is infinite, $\lim_{\lVert \widehat{\bm{\theta}} \rVert^2 \rightarrow \infty}  P(\bm{\theta} \le \bm{0} \mid \bm{y}) = 0,$ as required. \\

Now we turn to the case where $\widehat{\theta}_i \rightarrow -\infty$ for all $i$, in which case we assume that $w_i < 0$ for all $i$. Then, a Fr\'echet bound ensures that
$$
P(\bm{\theta} \le \bm{0} \mid \bm{y}) = P(\theta_1 \le 0, \theta_2 \le 0, \, ... \, , \theta_{r_1} \le 0 \mid \bm{y} ) \ge \sum_{i=1}^{r_1} P(\theta_i \le 0 \mid \bm{y} ) - (r_1-1).
$$
Therefore, 
$$
\lim_{\lVert \widehat{\bm{\theta}} \rVert^2 \rightarrow \infty} P( \theta_i \ge 0 \mid \bm{y}) = 0, \, 1 \le i \le r_1 \Rightarrow   \lim_{\lVert \widehat{\bm{\theta}} \rVert^2 \rightarrow \infty}  P(\bm{\theta} \le \bm{0} \mid \bm{y}) = 1.
$$
Then, we can work with the conditional probabilities exactly as we did for the previous case:
\begin{align*}
\lim_{\lVert \widehat{\bm{\theta}} \rVert^2 \rightarrow \infty}  P( \theta_i \ge 0 \mid \bm{y}) &=  \lim_{\lVert \widehat{\bm{\theta}} \rVert^2 \rightarrow \infty}  \int_0^\infty P(\theta_i \ge 0 \mid g, \bm{y}) p(g \mid \bm{y}) \, \mathrm{d} g \\
 &\lesssim\lim_{\lVert \widehat{\bm{\theta}} \rVert^2 \rightarrow \infty}  \frac{[\hatth'\Om^{-1} \hatth]^{-(n-p)/2}}{p(\bm{y})} \int_\varepsilon^\infty \frac{ P( T_{n-p} \ge -\bm{m}_j / \sqrt{\bm{S}_{jj}} ) \, (g+d_u)^{(n-p-r_1)/2}}{[(g+d_u) \mathsf{SSE}_1/\hatth'\Om^{-1} \hatth+ 1 ]^{(n-p)/2}}  \, \pi(\mathrm{d} g).
\end{align*}
Since $-\bm{m}_j$ is positive, the subsequent steps in the proof for the previous case allow us to conclude that $\lim_{\lVert \widehat{\bm{\theta}} \rVert^2 \rightarrow \infty}  P( \theta_i \ge 0 \mid \bm{y}) = 0,$ as required. \\
 
\noindent \underline{Necessary condition:}

In the sequel we assume that there is at least one $i$ such that $\widehat{\theta}_i \rightarrow + \infty$. The case where all coordinates go to $-\infty$ can be dealt with the same way we did for the sufficient condition. We can write:
\begin{align*}
 \lim_{\lVert \widehat{\bm{\theta}} \rVert^2 \rightarrow \infty} P( \bm{\theta} \le \bm{0} \mid \bm{y} ) &=  \lim_{\lVert \widehat{\bm{\theta}} \rVert^2 \rightarrow \infty}  \frac{[\hatth'\Om^{-1} \hatth]^{(n-p)/2} \int_0^\infty P(\th \le \bm{0} \mid g, \bm{y}) p(\bm{y} \mid g) \pi(\mathrm{d} g)}{[\hatth'\Om^{-1} \hatth]^{(n-p)/2} p(\bm{y}) }.
\end{align*}
First, we show that the limit of the numerator is bounded away from $0$. Applying Fatou's lemma and one of the bounds in Lemma~\ref{boundh},
\begin{align*}
 \lim_{\lVert \widehat{\bm{\theta}} \rVert^2 \rightarrow \infty} \frac{\int_0^\infty P(\th \le \bm{0} \mid g, \bm{y}) p(\bm{y} \mid g) \pi(\mathrm{d} g) }{[\hatth'\Om^{-1} \hatth]^{-(n-p)/2}}
 &\ge \int_{0}^\infty \lim \inf_{\lVert \widehat{\bm{\theta}} \rVert^2 \rightarrow \infty}  \frac{P(\th \le \bm{0} \mid g, \bm{y}) \,   h(g)} {[\hatth'\Om^{-1} \hatth]^{(n-p)/2}} \, \pi(\mathrm{d} g) \\
 &\gtrsim \int_{0}^\infty (g+d_l)^{(n-p-r_1)/2} \, \lim \inf_{\lVert \widehat{\bm{\theta}} \rVert^2 \rightarrow \infty}  \,  P(\th \le \bm{0} \mid g, \bm{y})  \pi(\mathrm{d} g),
\end{align*}
and for any $g$,
$$
\lim \inf_{\lVert \widehat{\bm{\theta}} \rVert^2 \rightarrow \infty}  \,  P(\th \le \bm{0} \mid g, \bm{y})  = \lim \inf_{\lVert \widehat{\bm{\theta}} \rVert^2 \rightarrow \infty}  \,  P(\bm{\xi} \le \bm{0} \mid g, \bm{y})
$$
where $\bm{\xi}$ is a multivariate Student-$t$ as in Lemma~\ref{inconBF2}. Lemma~\ref{inconBF2} shows that $ P(\bm{\xi} \le \bm{0} \mid g, \bm{y})$ is bounded away from 0, which implies that the numerator is bounded away from 0,
as claimed. A necessary condition for $\lim_{\normth^2 \rightarrow \infty} P(\bm{\theta} \le \bm{0} \mid y) = 0$ is that $\lim_{\normth^2 \rightarrow \infty} [\hatth'\Om^{-1} \hatth]^{(n-p)/2}  p(\bm{y}) = \infty$ which, as we saw in the proof of the sufficient condition, is equivalent to
$$
\int_0^\infty (g+1)^{(n-p-r_1)/2} \pi(\mathrm{d} g) = \infty,
$$
as required.

\section{Proof of Lemma \ref{inconBF3}}\label{proof.inconBF3}
The second part of the Bayes factor in \eqref{B10onesided2} can be expressed as $\left(P_{\pi}(\bm{\theta}\le\textbf{0} \mid \textbf{y})^{-1}-1\right)=\frac{k(\bm\Theta_1)}{k(\bm\Theta_0)}$, where
\[
k(\bm\Theta_t)=\int_{\boldsymbol{\theta}\in\bm\Theta_t} \left(\nu s^2+s_{\textbf{y}}^2+\left(\boldsymbol{\theta}-\hat{\boldsymbol{\theta}}\right)'\textbf{X}_{\bm\theta}'\boldsymbol{\Sigma}^{-1}\textbf{X}_{\bm\theta}\left(\boldsymbol{\theta}-\hat{\boldsymbol{\theta}}\right)\right)^{-\frac{n-r_2+\nu}{2}}
N_{\boldsymbol{\theta}}(\textbf{0},\bm\Omega,\bm\Theta_t)
d\boldsymbol{\theta},
\]
and $N_{\bm{\theta}}(\textbf{0},\bm\Omega,\bm\Theta_t)$ denotes a truncated multivariate normal density for $\bm\theta$ with mean $\textbf{0}$ and covariance matrix $\bm\Omega$, truncated in the subspace $\bm\Theta_t$ for $t=0$ or 1. Exactly as in the proof of Lemma~\ref{lemma.indconj} it can be shown that $k(\bm\Theta_t) = \left(\nu s^2+s_{\textbf{y}}^2+\hat{\boldsymbol{\theta}}'\textbf{X}_
{\bm\theta}'\boldsymbol{\Sigma}^{-1}\textbf{X}_{\bm\theta}\hat{\boldsymbol{\theta}}\right)^{-\frac{n-r_2+\nu}{2}}(1+o(1))$ in the limit, so that $\left(P_{\pi}(\bm{\theta}\le\textbf{0} \mid \textbf{y})^{-1}-1\right)\rightarrow 1$.

\section{Proof of Lemma \ref{lemmaNonnestedAdaptive_g}}
\label{proofNonnestedAdaptive_g}
The marginal posterior of $\bm\theta$ in the joint space has a multivariate Student t distribution with mean $\frac{g}{g+1}\hat{\bm{\theta}}$, scale matrix $(n-r_2)^{-1}(s_{\textbf{y}}^2+(g+1)^{-1}\hat{\bm{\theta}}'(\textbf{X}_{\bm\theta}'\bm{\Sigma}^{-1}\textbf{X}_{\bm\theta})\hat{\bm{\theta}})\frac{g}{g+1}
(\textbf{X}_{\bm\theta}'\bm{\Sigma}^{-1}\textbf{X}_{\bm\theta})^{-1}$, and $n-r_2$ degrees of freedom.
A change of variables to $\bm\xi=\frac{g+1}{g}$, results in a multivariate Student t distribution with mean $\hat{\bm{\theta}}$, scale matrix $(n-r_2)^{-1}((1+g^{-1})s_{\textbf{y}}^2+g^{-1}\hat{\bm{\theta}}'(\textbf{X}_{\bm\theta}'\bm{\Sigma}^{-1}\textbf{X}_{\bm\theta})\hat{\bm{\theta}})
(\textbf{X}_{\bm\theta}'\bm{\Sigma}^{-1}\textbf{X}_{\bm\theta})^{-1}$, and degrees of freedom $n-r_2$. Note that the posterior probability is invariant under this transformation, i.e., $P_{\pi}(\bm\theta\le\textbf{0}|\textbf{y})=P_{\pi}(\bm\xi\le\textbf{0}|\textbf{y})$. Furthermore it important to note that the factor $(1+g^{-1})s_{\textbf{y}}^2+g^{-1}\hat{\bm{\theta}}'(\textbf{X}_{\bm\theta}'\bm{\Sigma}^{-1}\textbf{X}_{\bm\theta})\hat{\bm{\theta}}$ in the scale matrix of $\bm\xi$ is a monotonically decreasing function of $g$. Now it is easy to see that if $\hat{\bm\theta}\le\textbf{0}$, $P_{\pi}(\bm\xi\le\textbf{0}|\textbf{y})$ monotonically increases as the scales decrease, and if $\hat{\bm\theta}\not\le\textbf{0}$, $P_{\pi}(\bm\xi\le\textbf{0}|\textbf{y})$ monotonically decreases as the scales decrease. Thus, in order to maximize $B_{01}$ if $\hat{\bm\theta}\le\textbf{0}$, and maximize $B_{10}$ if $\hat{\bm\theta}\not\le\textbf{0}$, we have to let $g$ go to $\infty$. For completeness note that the marginal posterior of $\bm\theta$ in the joint space with a multivariate Student t distribution with mean $\hat{\bm{\theta}}$, scale matrix $(n-r_2)^{-1}s_{\textbf{y}}^2(\textbf{X}_{\bm\theta}'\bm{\Sigma}^{-1}\textbf{X}_{\bm\theta})^{-1}$, and $n-r_2$ degrees of freedom, in the limit as $g\rightarrow\infty$. Thus, even though a (data based) adaptive prior is considered, the choice of $g$ that maximizes the Bayes factor does not depend on the data. Note that taking the limit $g\rightarrow\infty$ was already considered by \cite{Mulder:2014a} but not in the context of an adaptive prior.

\bibliographystyle{apacite}
\bibliography{refs_mulder}

\end{document}